\documentclass{amsart}

\usepackage[utf8]{inputenc}
\usepackage{amsmath,amssymb,amsthm,units, stmaryrd,stackrel,relsize,bm}
\usepackage{enumitem}
\usepackage{mathdots}
\usepackage{bussproofs}
\usepackage{tikz}
\usetikzlibrary{calc,arrows, decorations.pathmorphing, matrix,fit}
\usepackage{float}
\usepackage{setspace}

\newtheorem{theorem}{Theorem}

\newtheorem{lemma}[theorem]{Lemma}

\newtheorem{claim}[theorem]{Claim}
\newtheorem{corollary}[theorem]{Corollary}
\newtheorem{proposition}[theorem]{Proposition}
\newtheorem{question}[theorem]{Question}

\newtheorem{thmx}{Theorem}

\theoremstyle{definition}
\newtheorem{definition}[theorem]{Definition}

\theoremstyle{remark}
\newtheorem{remark}[theorem]{Remark}

\usepackage[
    colorlinks=true, 
    citecolor=red,
    linkcolor=blue,
    urlcolor=blue]{hyperref}

\newcommand\zfc{{\mathsf{ZFC}}}

\newcommand\kp{{\mathsf{KP}}}

\newcommand\Ord{{\mathsf{Ord}}}

\newcommand\dom{{\text{dom}}}

\newcommand\wo{{\textsc{wo}}}

\newcommand\ca{{\mathsf{CA_0}}}

\newcommand\aca{{\mathsf{ACA_0}}}
\newcommand\rca{{\mathsf{RCA_0}}}
\newcommand\eca{{\mathsf{ECA_0}}}

\makeatletter
\def\Ddots{\mathinner{\mkern1mu\raise\p@
\vbox{\kern7\p@\hbox{.}}\mkern2mu
\raise4\p@\hbox{.}\mkern2mu\raise7\p@\hbox{.}\mkern1mu}}
\makeatother

\begin{document}
\title{The $\Pi^1_2$ Consequences of a Theory}
\author{J. P. Aguilera}
\address{Institute of Discrete Mathematics and Geometry, Vienna University of Technology. Wiedner Hauptstra{\ss}e 8--10, 1040 Vienna, Austria \textit{and} Department of Mathematics, University of Ghent. Krijgslaan 281-S8, B9000 Ghent, Belgium}
\email{aguilera@logic.at}
\author{F. Pakhomov}
\address{Department of Mathematics, University of Ghent. Krijgslaan 281-S8, B9000 Ghent, Belgium \textit{and} Steklov Mathematical Institute of the Russian Academy of Sciences. Ulitsa Gubkina 8, Moscow 119991, Russia.}
\email{fedor.pakhomov@ugent.be}
\thanks{The authors would like to thank Anton Freund, Harvey Friedman, Michael Rathjen, and Andreas Weiermann for fruitful conversations related to the subject of the article.
The first author was partially supported by FWF grant I4513N and FWO grant 3E017319. The second author was partially supported by FWO grant G0F8421N}

\begin{abstract}
We develop the abstract framework for a proof-theoretic analysis of theories with scope beyond ordinal numbers, resulting in an analog of \mbox{Ordinal Analysis} aimed at the study of theorems of complexity $\Pi^1_2$. This is done by replacing the use of ordinal numbers by particularly uniform, wellfoundedness preserving functors in the category of linear orders.

Generalizing the notion of a \emph{proof-theoretic ordinal},
we define the functorial $\Pi^1_2$ norm of a theory and prove its existence and uniqueness for $\Pi^1_2$-sound theories.
From this, we further abstract a definition of the \emph{$\Sigma^1_2$- and $\Pi^1_2$-soundness ordinals of a theory}; these quantify, respectively, the maximum strength of true $\Sigma^1_2$ theorems and minimum strength of false $\Pi^1_2$ theorems of a given theory. We study these ordinals, developing a proof-theoretic classification theory for recursively enumerable extensions of $\aca$

Using techniques from infinitary and categorical proof theory, generalized recursion theory, constructibility, and forcing, we prove that an admissible ordinal is the $\Pi^1_2$-soundness ordinal of some recursively enumerable extension of $\aca$ if and only if it is not parameter-free $\Sigma^1_1$-reflecting. 
We show that the $\Sigma^1_2$-soundness ordinal of $\aca$ is $\omega_1^{ck}$ and characterize the $\Sigma^1_2$-soundness ordinals of recursively enumerable, $\Sigma^1_2$-sound extensions of $\Pi^1_1{-}\ca$.
\end{abstract}
\date{\today $\,$ (compiled)}
\clearpage
\keywords{Proof Theory, $\Pi^1_2$ Sentence, Ordinal Analysis, Ordinal Number, Dilator, Wellordering Principle}
\subjclass[2020]{03F15, 03F25 (primary), 18A99 03D20, 03D60, 03D65, 03E45 (secondary)}
\maketitle
\numberwithin{equation}{section}
\setcounter{tocdepth}{1}
\tableofcontents

\section{Introduction}
The purpose of this article is to develop an abstract framework to justify the prospect of a $\Pi^1_2$ proof-theoretic analysis of mathematical theories, extending the theory of Ordinal Analysis to objects of higher complexity.

The study of proof-theoretic analysis of theories goes back to Gentzen \cite{Ge36} who proved the consistency of Peano Arithmetic using an argument involving induction applied to a computable process along an ordering  of length
\[\varepsilon_0 = \sup\{\omega,\omega^\omega, \hdots\},\]
where $\omega$ denotes the order-type of $\mathbb{N}$. This result provides further context for G\"odel's \cite{Goe31} impossibility proof for Hilbert's second problem and indeed provides an answer conditioned on a quasi-finitistic component.
Moreover, this type of proof-theoretic analysis of a theory yields insight into the set of its consequences and has many applications, e.g., in relation to unprovability and combinatorics. For some examples, we refer the reader to the work of Kirby-Paris \cite{KP82}, Paris-Harrington \cite{PH77}, or Ketonen-Solovay \cite{KS81}.

Takeuti \cite{Ta67} famously proved an analogous result for the subsystem $\Pi^1_1$ of analysis.
With time, the ideas underlying Gentzen's and Takeuti's proof solidified into the field of \emph{ordinal analysis}. General references are Pohlers \cite{Po08}, Sch\"utte \cite{Sch77}. The goal of ordinal analysis is to study a mathematical theory $T$ by isolating its finitary and infinitary components. All known ordinal analyses of theories that have been carried out additionally provide deep insight into the structural behaviour of the axioms of the theory and how they could be unravelled into potential direct proofs of a contradiction. 
We refer the reader to Rathjen \cite{Ra99a} for an overview.

The true content of these results, however, is difficult to state succinctly. Hence, it has become tradition in proof theory to speak of the \emph{proxy} problem of carrying out a $\Pi^1_1$-analysis of a theory $T$. Let $\wo(a)$ be the formula expressing that $a$ is a wellordering with field a subset of $\mathbb{N}$. Define
\begin{align*}
|T|_{\Pi^1_1} = \sup\{\alpha\in\Ord: \text{$T\vdash \wo(a)$ for some recursive linear order $a\cong\alpha$}\}.
\end{align*}
$|T|_{\Pi^1_1}$ is always defined if $T$ is a recursively enumerable, $\Pi^1_1$-sound theory. In this case, we additionally have $|T|_{\Pi^1_1}<\omega_1^{ck}$, i.e., $|T|_{\Pi^1_1}$ is a recursive ordinal. This is an immediate consequence of the $\Sigma^1_1$-boundedness theorem.  By convention, one usually defines $|T|_{\Pi^1_1} = \omega_1^{ck}$ if $T$ is $\Pi^1_1$-unsound (but recursively enumerable).

The problem becomes then to compute $|T|_{\Pi^1_1}$ for a specific theory $T$. This involves exhibiting an explicit recursive ordering of $\mathbb{N}$ of order-type $|T|_{\Pi^1_1}$. The proof that this fact holds will ideally involve a direct reduction of every provable wellorder of $T$ into the ordering exhibited. Since wellfoundedness is  a complete $\Pi^1_1$-property, this leads to a characterization of the $\Pi^1_1$ consequences of $T$. 

Historically, Proof Theory has encountered difficulties extending its scope of studies from $\Pi^1_1$ sentences to higher complexities, and this is perhaps in part due many of the techniques relying on the use of ordinals, whose extremely simple structure is both a blessing and a curse. Recently, however, attempts at replacing the use of ordinals with more complex categorical constructions have gained momentum and resulted in higher order counterparts of classical theorems (see e.g. Freund \cite{Fr19}, a categorical analog of the classical well-ordering theorems of Girard (unpublished), Friedman (unpublished), Marcone-Montalb\'an \cite{MaMo11}, Rathjen-Weiermann \cite{RaWe11}, and \cite{Agetal}, a categorical analog of the Kirby-Paris theorem).

The theory of $|T|_{\Pi^1_1}$ is an abstract framework which makes possible the prospect of $\Pi^1_1$-analyses of theories $T$. Among other things, it relies crucially on the existence of $|T|_{\Pi^1_1}$. In this article, we generalize these results to the class of $\Pi^1_2$ formulas. In particular, we give a robust definition of $|T|_{\Pi^1_2}$. 

\subsection{The $\Pi^1_2$ consequences of a theory}
In order to study $\Pi^1_2$ consequences of a theory, we need to abandon the idea of relying uniquely on ordinals, as the class of ordinal notations is $\Pi^1_1$. Instead, we need to work with dilators, which were first introduced by Girard \cite{Gi81}. We consider $\Ord$ as the category of wellorders where morphisms are strictly increasing functions. A \emph{dilator} is a functor on $\Ord$ which commutes with pullbacks and direct limits.
DIL can be regarded as a functor category, with natural transformations as morphisms.

\begin{definition}\label{DefinitionPi12Dilator}
Suppose $T$ is a theory, then $|T|_{\Pi^1_2}$ is the unique dilator $D^*$ up to bi-embeddability with the following properties, if it exists:
\begin{enumerate}
\item Suppose $T$ proves that $D$ is a recursive dilator, then $D$ embeds into $D^*$; and
\item suppose $\hat D$ satisfies (1), then $D^*$ embeds into $\hat D$, and moreover the diagram commutes.
\end{enumerate} 
\end{definition}

\begin{figure}[h!]
\begin{tikzpicture}
  \matrix (m) [matrix of math nodes,row sep=3em,column sep=4em,minimum width=2em]
  {
     D_0 &   & \\
      \vdots   & {|T|_{\Pi^1_2}}  & \hat D\\
     D_i & &  \\};
  \path[-stealth]
    (m-1-1) edge (m-2-2)
    (m-3-1) edge (m-2-2)
    (m-1-1) edge [bend left] (m-2-3)
    (m-3-1) edge [bend right](m-2-3)
     (m-2-2) edge [densely dotted] (m-2-3);
\end{tikzpicture}
\caption{The universal property for $|T|_{\Pi^1_2}$. Here, the functors $D_j$ are the provable dilators of $T$ and the arrows represent natural transformations.}
\end{figure}
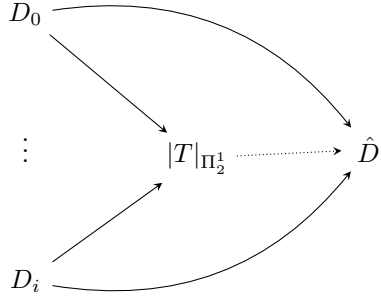

A dilator $D$ satisfying the first condition in Definition \ref{DefinitionPi12Dilator} is easy to find. In fact, this is immediate from the $\Sigma^1_1$-boundedness theorem for dilators, due independently to Girard-Normann \cite{GN85} and Kechris-Woodin \cite{KW91}. However, such a dilator will not necessarily satisfy the second condition. Since DIL is not wellordered by embeddability (unlike $\Ord$), it is not immediately clear that $|T|_{\Pi^1_2}$ can always exist. Our first theorem states that it does.

\begin{thmx}\label{thmxA}
Suppose that $T$ is a $\Pi^1_2$-sound, recursively enumerable extension of $\aca$. Then, $|T|_{\Pi^1_2}$ exists and is recursive.
\end{thmx}
The theorem is proved in \S  \ref{SectionNorm}, where we also derive some consequences of the definition and existence. Catlow \cite{Ca94} has previously carried out a study of the provable dilators of $\aca$. He finds a function on $\Ord$ which extensionally bounds all the provable dilators of $\aca$, i.e., it bounds the restriction of all such dilators to $\Ord$.
A similar line of research was pursued in \cite{PW21}, where the notion of a \emph{proof-theoretic dilator} was defined. This definition was again extensional.

In contrast, the definition of $|T|_{\Pi^1_2}$ is \emph{intensional}, i.e., functorial. This is crucial in order for it to relate back to the $\Pi^1_2$ consequences of $T$. Moreover, all the structure theory for dilators can be applied to it. In particular, it can be coded by a recursive subset of $\mathbb{N}$ in such a way that its extension to any ordinal $\alpha$ of any cardinality is uniquely determined. 

The computation of $|\aca|_{\Pi^1_2}$ is somewhat more involved and will be reported in a forthcoming article \cite{AgPa}. Nonetheless, even without knowing what $|T|_{\Pi^1_2}$ is for a given $T$, we can deduce a fair amount of information about and from it. In particular, we prove the following extensional description which ties the theory developed here with the work of \cite{Ca94} and \cite{PW21}.
\begin{thmx}
Suppose $T$ is $\Pi^1_2$-sound and let $D = |T|_{\Pi^1_2}$. Then, for every recursive wellordering $\alpha$ of $\mathbb{N}$, we have 
\[D(\alpha) = |T + \wo(\alpha)|_{\Pi^1_1}.\]
\end{thmx}

\subsection{A proof-theoretic classification theory}
The definition of $|T|_{\Pi^1_2}$ applies only to $\Pi^1_2$-sound theories. $\Pi^1_2$-unsound theories are of potential interest and appear e.g., in Reverse Mathematics, where they are used in gauging the strength of many mathematical theorems, such as Martin's Borel Determinacy \cite{Ma75}. Friedman \cite{Fr71} showed that $\Sigma^0_5$-determinacy is not provable in Second-Order Arithmetic, Martin (unpublished) improved this result to $\Sigma^0_4$, and Montalb\'an-Shore \cite{MoSh11} to $\bigcup_{n\in\mathbb{N}} n{-}\Pi^0_3$. Proofs of this type of result generally begin with the assumption ``Suppose there is no transitive model of $T$.''

Partly motivated by this,
in \S  \ref{SectUnsoundness}, we study the natural attempt to extend $|T|_{\Pi^1_2}$ to $\Pi^1_2$-unsound theories. Of course, no such attempt can succeed. 
However, for a specific theory $T$, we can keep record of how far the attempt goes before breaking down and 
synthesize this value into an
ordinal measure $o^1_2(T)$ which quantifies how close $T$ is to being a $\Pi^1_2$-sound theory.  This ordinal is defined in \S  \ref{SectUnsoundness} and it is the least ordinal $\alpha$ such that $D(\alpha)$ is illfounded for some provable dilator $D$ of $T$.

The ordinal $o^1_2(T)$ can be used to obtain useful information about $T$. For instance, it allows us to separate all recursively enumerable extensions of $\aca$ into four categories, in a sort of proof-theoretic counterpart to the classification theory which has proven to be extremely fruitful in Model Theory.
Category A is comprised of theories $T$ with $o^1_2(T) = 0$; Category B is comprised of theories with $o^1_2(T)$ nonzero, but recursive; Category C is comprised of theories with $o^1_2(T)$ non-recursive; and Category D is comprised of theories with $o^1_2(T)$ undefined (in which case we write $o^1_2(T) = \infty$. Below, let $\mathsf{Bool}(\Pi^1_1)$ denote the class of all Boolean combinations of $\Pi^1_1$ sentences.

\begin{figure}[h!]
\begin{tikzpicture}
    \fill[fill=red!20] (-5.05,0.15) rectangle (-4.9,-0.15);
    \fill[fill=orange!20] (-4.05,0.15) rectangle (-2,-0.15);
    \fill[fill=yellow!40] (-2,0.15) rectangle (2.05,-0.15);
    \fill[fill=green!30] (4.9,0.15) rectangle (5.05,-0.15);

  \begin{scope}[thick]
    \draw[|-)]                 (-5,0) node[left]{}  -- (-4,0) 
                                        node[right]{};
    \draw[-)]                 (-4,0) node[left]{}  -- (-2,0) 
                                        node[right]{};
    \draw[-)]                 (-2,0) node[left]{}  -- (2,0) 
                                        node[right]{};
    \draw[-)]                 (2,0) node[left]{}  -- (5,0) 
                                        node[right]{};
     \draw[]             (-5,-0.3) node {$0$};
     \draw[]             (-4,-0.3) node {$\varepsilon_0$};
     \draw[]             (-2,-0.3) node {$\omega_1^{ck}$};
     \draw[]             (2,-0.3) node {$\delta^1_2$};
     \draw[]             (5,-0.3) node {$\infty$};
     \draw[]             (-5,0.4) node {A};
     \draw[]             (-3,0.4) node {B};
     \draw[]             (0,0.4) node {C};
     \draw[]             (5,0.4) node {D};
  \end{scope}

\end{tikzpicture}
\caption{The four categories of recursively enumerable extensions of $\aca$ according to their degree of $\Pi^1_2$-soundness.}
\end{figure}
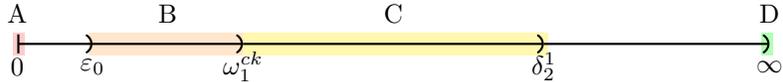
\clearpage

\begin{thmx}\label{thmxB}
Suppose that $T$ is a recursively enumerable extension of $\aca$. Then,
\begin{enumerate}
\item[A.] \label{thmxBA}The following are equivalent:
\begin{enumerate}
\item $T$ is in Category A, i.e., $o^1_2(T) = 0$;
\item $o^1_2(T) < |T|_{\Pi^1_1}$;
\item $T$ is not $\Pi^1_1$-sound.
\end{enumerate}
\item [B.] \label{thmxBB} The following are equivalent:
\begin{enumerate}
\item $T$ is in Category B, i.e., $o^1_2(T)$ is nonzero, but recursive;
\item $|T|_{\Pi^1_1} \leq o^1_2(T) < \omega_1^{ck}$;
\item $T$ is $\Pi^1_1$-sound, but $\mathsf{Bool}(\Pi^1_1)$-unsound.
\end{enumerate}
\item [C.] \label{thmxBC} The following are equivalent:
\begin{enumerate}
\item $T$ is in Category C, i.e., $o^1_2(T)$ is non-recursive;
\item $\omega_1^{ck}\leq o^1_2(T) < \delta^1_2$;
\item $T$ is $\mathsf{Bool}(\Pi^1_1)$-sound, but $\Pi^1_2$-unsound.
\end{enumerate}
\item[D.] \label{thmxBD} The following are equivalent:
\begin{enumerate}
\item $T$ is in Category D, i.e., $o^1_2(T) = \infty$;
\item $\delta^1_2 \leq o^1_2(T)$;
\item $T$ is $\Pi^1_2$-sound.
\end{enumerate}
\end{enumerate}
\end{thmx}

\subsection{The $\Pi^1_2$-Spectrum Problem}
Theorem \ref{thmxB} splits the recursively enumerable extensions $T$ of $\aca$ into four categories according to the value $o^1_2(T)$. Additionally, it shows that some values of $o^1_2(T)$ are impossible. For instance, it shows that
\[0 < o^1_2(T)\to |T|_{\Pi^1_1} < o^1_2(T),\]
so no theory can have value strictly between $0$ and $\varepsilon_0$. Similarly, it shows that $o^1_2(T)$ is always countable. 
It does not say e.g., whether $o^1_2(T)$ could be have arbitrarily large recursive values, or whether it could be admissible.

These questions lead to the \emph{Spectrum Problem}, which consists of identifying the possible ordinals which are of the form $o^1_2(T)$ for some recursively enumerable extension $T$ of $\aca$. The Spectrum Problem can be split into two subproblems, according as one is dealing with recursive or non-recursive ordinals, and thus with theories in Category B or C.
We solve the problem for theories in Category B in \S  \ref{SectRecursive}:
\begin{thmx}\label{thmxC}
Let $\alpha$ be a recursive ordinal. The following are equivalent:
\begin{enumerate}
\item $\alpha = o^1_2(T)$ for some recursively enumerable extension of $\aca$;
\item $\alpha = \varepsilon_{\gamma}$ for some recursive $\gamma$.
\end{enumerate}
\end{thmx}

Using the machineries of $\beta$-proofs and pointed Sacks forcing, we give an answer to the Spectrum Problem theories in Category C which solves the problem for admissible ordinals. This is done in \S  \ref{SectAdmissible}:
\begin{thmx}\label{thmxD}
Let $\alpha$ be an admissible ordinal or a limit of admissible  ordinals. The following are equivalent:
\begin{enumerate}
\item $\alpha = o^1_2(T)$ for some recursively enumerable extension of $\aca$;
\item $\alpha$ does not reflect every parameter-free $\Sigma^1_1$ sentence.
\end{enumerate}
\end{thmx}

\subsection{The $\Sigma^1_2$-soundness ordinal}
In \S \ref{SectSigma}, we study a sort of ``positive'' dual to the $\Pi^1_2$-soundness ordinal, denoted $s^1_2(T)$. While $o^1_2(T)$ measures the ``truth complexity'' of the \emph{false} $\Pi^1_2$ statements provable by $T$, $s^1_2(T)$ measures the complexity of the \emph{true} $\Sigma^1_2$ sentences provable by $T$. While this ordinal has higher intrinsic motivation, it does not seem to allow for a proof-theoretic classification theory like $o^1_2(T)$ does.
We prove various results about $s^1_2(T)$, the most important of which we summarize here:

\begin{thmx}\label{thmxE}
Let $T$ be a $\Sigma^1_2$-sound theory.
\begin{enumerate}
\item $s^1_2(T) <\delta^1_2$ and the bound is optimal;
\item $s^1_2(T) \leq o^1_2(T)$, and both equality and inequality are possible;
\item $s^1_2(\aca) = s^1_2(\kp) = \omega_1^{ck}$;
\item  $s^1_2(\Pi^1_1$-$CA_0) = \omega_\omega^{ck}$;
\item $s^1_2(\Pi^1_2$-$CA_0)$ is the least ordinal stable to the least non-projectible ordinal.
\end{enumerate}
\end{thmx}

We also give a characterization of the ordinals of the form $s^1_2(T)$ for some recursively enumerable, $\Sigma^1_2$-sound extension of $\Pi^1_1$-CA$_0$.

\section{Preliminaries}\label{SectionPreliminaries}
All our notation and definitions are standard.
Nonetheless, in this section we collect a list of preliminary  definitions and known results.
The most important concepts from Proof Theory which we will need were defined in the introduction, but we will require some notions from various areas, which we recall below. We will, however, speak of proofs. We will think of proofs as carried out in a sequent calculus, e.g., as in Takeuti \cite{Ta75}, but the choice of formalism will be inconsequential.\\

\paragraph{\textbf{Subsystems of Second-Order Arithmetic}}
We will deal with subsystems of Second-Order Arithmetic. For background, we refer the reader to Simpson \cite{Simpson}. The main system we will deal with is $\aca$, Arithmetical Comprehension. This is the second-order analog of Peano Arithmetic containing the Induction axiom, as well as the schema asserting the existence of every first-order definable subset of $\mathbb{N}$ (with parameters). We will also consider other subsystems of Second-Order Arithmetic, including $\Pi^1_1${-}$\ca$, $\Pi^1_2${-}$\ca$. Moreover, we will consider theories in the language of set theory, such as $\kp$, Kripke-Platek set theory. We will sometimes be unclear about whether we deal with theories in the language of Second-Order Arithmetic, or in the language of set theory, since the two languages can easily be translated into one another and this should lead to no confusion.  Nonetheless, the theories $\Pi^1_1${-}$\ca$, $\Pi^1_2${-}$\ca$ always refer to comprehension for sets of natural numbers (and never to any kind of second-order comprehension in the language of set theory).

We will speak of countable objects (such as functions on the natural numbers or models in the language of set theory) which can be coded by sets of natural numbers. In order to make the exposition simpler, we will often identify objects with their codes. We will also identify sets of natural numbers with real numbers.

We will make use of the notion of an $\omega$-model, which is a model in the language of Second-Order Arithmetic (or set theory) whose natural numbers are isomorphic to $\mathbb{N}$. \\

\paragraph{\textbf{Recursion theory}} 
We will make use of basic Admissible Recursion Theory. The standard reference is Barwise \cite{Ba75}.
A transitive set is \emph{admissible} if it is a model of $\kp$. In particular, we will consider admissible sets of the form $L_\alpha$, where $L$ is G\"odel's constructible universe. Ordinals $\alpha$ such that $L_\alpha$ are admissible are called \emph{admissible ordinals}. 
An ordinal is \emph{recursive} if it is the order-type of a recursive wellordering of $\mathbb{N}$. The least non-recursive ordinal is called $\omega_1^{ck}$. By abuse of notation, we will often identify ordinals with ordinal notations for them, but only if this causes no confusion. We generally employ greek letters $\alpha,\gamma,\eta,\xi,\zeta$ for ordinal numbers or their codes. When dealing with an ordinal code which might possibly be illfounded, or if the distinction between an ordinal and its code is relevant, we employ letters from the Roman alphabet, mainly $a$. If $\alpha$ is an ordinal, we use $\alpha^+$ to denote the smallest admissible ordinal greater than $\alpha$.

Given a real number $x$ (a set of natural numbers), we define $L[x]$, the constructible hierarchy relativized to $x$, as usual. We denote by $\omega_1^x$ the least ordinal which is not recursive relative to $x$. For each $x$, $\omega_1^x$ is the least $\alpha$ such that $L_{\alpha}[x]$ is admissible. We say that $x$ is \emph{hyperarithmetical} in $y$ ($y\leq_{hyp} x$) if $y\in L_{\omega_1^x}[x]$. For more on Generalized Recursion Theory, we refer the reader to Sacks \cite{Sa90}. For more on constructibility, we refer the reader to Jech \cite{Je03} or Jensen \cite{Je72}.

We will deal extensively with trees $\mathbb{N}$, viewed as finite sequences of natural numbers, with infinite branches corresponding to infinite strings of natural numbers. A tree is \emph{wellfounded} if it has no infinite branches. The Kleene-Brouwer ordering on elements of a tree is defined by setting $s < t$ if $s$ is a proper extension of $t$ or, if letting $i$ be least such that $s_i \neq t_i$, we have $s_i < t_i$. A tree is wellfounded if and only if its Kleene-Brouwer ordering is wellfounded.  \\

\paragraph{\textbf{Dilators}}
We will make extensive use of the basic theory of dilators. For more background, we refer the reader to Girard \cite{Gi81} or Girard-Normann \cite{GN85}.
As mentioned in the introduction, we consider the category $\wo$ of wellorders with strictly increasing embeddings as morphisms. By abuse of notation, we will often identify $\Ord$ and $\wo$.
A dilator is a functor in this category which preserves pullbacks and direct limits. 
DIL is the functor category of dilators, with morphisms as natural transformations. We sometimes call natural transformations \emph{embeddings}.

A dilator $D$ is \emph{countable} if it maps countable wellorders to countable wellorders.
By commutation with direct limits, dilators are uniquely (up to natural isomorphism) determined by their action on finite linear orders (every linear order is the direct limit of its finite suborders).
Thus, countable dilators can be coded by real numbers and thus it makes sense to talk about them in the context of Second-Order Arithmetic. 

There are several ways to do this. For the sake of definiteness, we can think of coded dilators as the restriction of its domain to the category of finite ordinals and strictly increasing maps between them. That is a coded dilator is a family of well-orders $\langle D(n)\mid n\in \mathbb{N}\rangle$ and maps between well orders $\langle D(f)\colon D(n)\to D(m)\mid n,m\in\mathbb{N}\text{ and }f\colon n\to m\text{ is strictly increasing}\rangle$. We fix some construction of orders $D(A)$ for all orders $A$.

A pre-dilator is a functor on the category of linear orders which commutes with direct limits and pullbacks. Hence, a dilator is a pre-dilator which preserves wellfoundedness.
The set of all codes of dilators is $\Pi^1_2$-complete, but the set of codes of pre-dilators is recursive.  The fact that the set of all dilators is $\Pi^1_2$-complete can be proved by a tree-construction similar to the proof of Shoenfield absoluteness, using Kleene-Brouwer orderings. 
Hence, this fact is provable in $\aca$ (this result is due to Catlow \cite{Ca94}).
A coded dilator $D$ is \emph{recursive} if its code is recursive; equivalently, if the functions
\[f\mapsto D(f)\]
and
\[n\mapsto D(n),\]
for finite $n$ and strictly increasing $f: n\to m$ between naturals, are recursive. 
We will carry out the usual abuse of notation and identify codes of dilators with the dilators themselves, if this causes no confusion. In particular, we will henceforth denote the image of (a code of) an ordinal or a function under $D$ by $D(\alpha)$ or $D(f)$, and forget about the notation $D_i$ employed earlier.

Dilators can also be represented as denotation systems for ordinal numbers. These are formed by a collection of terms $t(x_1,\hdots, x_n)$, each with a fixed (possibly null) arity, together with comparison rules for terms. The comparison rules must specify which of  $t(x_1,\hdots, x_n)$ and $s(y_1,\hdots, y_m)$ is bigger, whenever we are given:
\begin{enumerate}
\item an $n$-ary term $t$, 
\item an $m$-ary term $s$,
\item ordinals $x_1 < \dots < x_n$,
and
\item ordinals $y_1 < \dots < y_n$.
\end{enumerate}  
The rules must not depend on the individual ordinals chosen, but only on their relative orderings.
We will make use of both ways of thinking about dilators. For a proof of the equivalence between the two definitions, we refer the reader to Girard \cite{Gi81}.
\\

\paragraph{\textbf{$\beta$-logic}}
$\beta$-logic takes place in the language of first-order logic with an added relation symbol $<$\footnote{It is common to formulate $\beta$-logic in the setting where the sort of ordinals $o$ is just one of the sorts, which is natural for some of the applications of the $\beta$-logic. However in the present paper it will be sufficient to consider the one-sorted variant of $\beta$-logic}. 
A $\beta$-structure is a model in this language where $<$ is interpreted as a wellordered relation. A sentence is $\beta$-satisfiable if it has such a model and $\beta$-valid if it has no such countermodel.
There is a corresponding notion of $\beta$-proof for $\beta$-structures, and this is also functorial. 

Suppose $\alpha$ is an ordinal. An $\alpha$-proof is a proof $P(\alpha)$ in first-order logic except that
\begin{enumerate}
\item no variables appear in $P(\alpha)$,
\item the language has constants $c_\iota$ for each $\iota<\alpha$, and axioms $c_\iota < c_\xi$ whenever $\iota<\xi$, and $c_\iota \leq c_\xi$ whenever $\iota
\leq\xi$.
\item $P(\alpha)$ is allowed to use the infinitary \emph{$\alpha$-rule}: from $A(c_\iota)$ for each $\iota<\alpha$, conclude $\forall x\, A(x)$.
\end{enumerate}
Given an embedding $f:\alpha\to\alpha'$, an $\alpha$-proof $P(\alpha)$ and an $\alpha'$-proof $P(\alpha')$, one can attempt to define an embedding $P(f):P(\alpha)\to P(\alpha')$ from sequents in  $P(\alpha)$ to sequents in $P(\alpha')$ such that 
\begin{enumerate}
\item $P(f)$ preserves the predecessor relation between sequents,
\item $P(f)(\Gamma)$ is a sequent identical to $\Gamma$, except that each constant $c_\iota$ has been replaced by $c_{f(\iota)}$,
\item $P(f)$ maps the conclusion of $P(\alpha)$ to the conclusion of $P(\alpha')$.
\end{enumerate} 
Note that these embeddings are functorial (they respect identity and composition). 

A $\beta$-proof is a family  $\{P(\alpha):\alpha\in\Ord\}$ such that $P(\alpha)$ is an $\alpha$-proof for each $\alpha$ and
every $f: \alpha\to\alpha'$  extends to an embedding $P(f)$ as above. If so, the system 
\[\Big\{\{P(\alpha):\alpha\in\Ord\}, \{P(f): \text{$f$ is a morphism in $\Ord$}\}\Big\}\]
commutes with direct limits and pullbacks. 
Thus, a $\beta$-proof is an analog of dilators for the category of infinitary proofs, though it is not necessary for us to explicitly consider this category.

Like a dilator, a $\beta$-proof is uniquely determined by the family $\{P(n): n\in\mathbb{N}\}$. 
A \emph{$\beta$-pre-proof} is a family $P = \{P(n): n\in\mathbb{N}\}$ which admits embeddings $P(f)$ whenever $f$ is an embedding between natural numbers. 
If so, then $P$ can be uniquely extended to infinite arguments $\alpha$ in such a way that the family $\{P(\alpha):\alpha\in\Ord\}$ still admits embeddings $P(f)$ whenever $f$ is an embedding in $\Ord$. This family will satisfy commutation with direct limits and pullbacks. 
It need not, however, satisfy that $P(\alpha)$ is a wellfounded proof tree for each $\alpha\in\Ord$.

If a $\beta$-pre-proof $P$ as above has the additional property that $P(\alpha)$   is wellfounded for each $\alpha$, then we identify it with the corresponding $\beta$-proof and with codes for it. We generally denote (codes for) $\beta$-proofs and $\beta$-pre-proofs by the letter $P$.

We will need the completeness theorem for $\beta$-logic, due to Girard. Since we will need a specific form of it, we sketch the proof.
\begin{theorem}[Girard] \label{TheoremBetaCompleteness}
Let $\varphi$ be a closed sentence of first-order logic, with a distinguished binary relation $<$. Then over $\aca$, we can effectively find a cut-free $\beta$-pre-proof $P$ such that the following are equivalent for each $\gamma$:
\begin{enumerate}
\item $\varphi$ is valid in all $\beta$-structures in which $<$ is interpreted as membership in $\gamma$,
\item $P(\gamma)$ is a wellfounded proof tree.
\end{enumerate}
In particular, $\varphi$ is $\beta$-valid if and only if $P$ is a $\beta$-proof.
\end{theorem}
\proof[Proof Sketch]
We reason in $\aca$. We need to begin by verifying soundness for cut-free $\beta$-proofs. Fix a cut-free $\beta$-proof $P$ and an ordinal $\gamma$.
For each $\beta$-structure $M$ as in the statement of the theorem, we check by induction on the tree $P(\gamma)$ that $M$ satisfies every sequent in $P(\gamma)$. Since $P(\gamma)$ is cut-free, all formulas in $P(\gamma)$ have complexity bounded by that of $\varphi$. Using a partial truth predicate for $M$ of sufficient complexity, define the set $X$ of all sequents $\Gamma$ in $P(\gamma)$ such that $M\not\models\Gamma$. We claim that there is $\Gamma\in X$ which is maximal with respect to the tree ordering of $P(\gamma)$. Otherwise, every element in $X$ has a successor, so $X$ is a tree with no terminal nodes. Since $P(\gamma)$ comes from a $\beta$-pre-proof, it is recursively bounded, so by K\"onig's lemma there is an infinite branch through $X$, contradicting wellfoundedness. 

We have shown that there is $\Gamma\in X$ which is maximal with respect to the tree ordering of $P(\gamma)$. But for this $\Gamma$, we have $M\not\models\Gamma$, yet $M$ satisfies every premise of $\Gamma$, which is impossible.

In order to prove completeness, fix a formula $\varphi$ and an ordinal $\gamma$. We consider the (possibly illfounded) cut-free proof tree obtained via the usual Sch\"utte-type completeness proof for $\omega$-logic. Such a tree is arithmetical in $\gamma$, so $\aca$ is enough for its existence.

If the tree has an infinite branch $b$, then one can use $\aca$ to construct a countermodel from it as in the proof of completeness for $\omega$-logic. The assignment of values to atomic formulas, as well as the proof that the structure obtained is indeed a countermodel, requires searching through the nodes in $b$, so the structure is arithmetical in $b$ and $\gamma$.

If the tree has no infinite branches, then it is a $\gamma$-proof, so by soundness it holds in every $\beta$-structure in which $<$ is interpreted as membership in $\gamma$.
This construction is functorial, so this defines a $\beta$-pre-proof, as desired.
\endproof

\section{The $\Pi^1_2$-norm of a theory} \label{SectionNorm}
The subject of study in this section will be $\Pi^1_2$-sound theories.
We will define the $\Pi^1_2$-norm of a theory $|T|_{\Pi^1_2}$ and prove essential facts about it, including its existence. 

\begin{definition} \label{DefPi12}
Suppose $T$ is a theory, then $|T|_{\Pi^1_2}$ is defined as the unique dilator $D^*$ up to bi-embeddability with the following properties, if it exists:
\begin{enumerate}
\item \label{DefPi12Cond1} Suppose $T$ proves that $D$ is a recursive dilator, then $\bar D$ embeds into $D^*$; and
\item \label{DefPi12Cond2} suppose $\hat D$ satisfies \eqref{DefPi12Cond1}, then $D^*$ embeds into $\hat D$ and the embeddings commute.
\end{enumerate} 
\end{definition}

\begin{theorem} \label{TheoremPi12Dilator}
Suppose $T$ is $\Pi^1_2$-sound and extends $\eca$. Then, $|T|_{\Pi^1_2}$ is defined. Moreover, if $T$ is recursively enumerable, then $|T|_{\Pi^1_2}$ is recursive.
\end{theorem}
\proof
Let $\mathcal{D} = \{D_0, D_1, D_2, \hdots\}$ be an enumeration of all recursive pre-dilators $D$ such that
\[T\vdash \text{``$D$ is a dilator.''}\]
By $\Pi^1_2$-soundness, each $D_i$ is a dilator. We may form the dilator
\[D^* = \sum_{i\in\mathbb{N}} D_i\]
by setting
\begin{align*}
D^*(\alpha) &= \sum_{i\in\mathbb{N}} D_i(\alpha)
\end{align*}
and, given $f: n\to m$ increasing,
\begin{align*}
D^*(f) : \sum_{i\in\mathbb{N}} D_i(n) &\to \sum_{i\in\mathbb{N}} D_i(m)\\
 \sum_{k<k^*}D_k(n) + l&\mapsto \sum_{k<k^*}D_k(m) + D_{k^*}(f)(l), &\text{ for $l < D_{k^*}(n)$.} 
\end{align*}

Clearly, $D^*$ satisfies condition \eqref{DefPi12Cond1} in Definition \ref{DefPi12}. 
We must show that it satisfies \eqref{DefPi12Cond2} as well. 

Thus, let $\hat D$ be such that each $D_i$ embeds into $\hat D$. We describe an embedding of $D^*$ into $\hat D$. For this, we make use of the perfect decomposition theorem for dilators of Girard \cite[Theorem 3.1.5]{Gi81}. It states that for every dilator $F$ there is a unique ordinal $\alpha$ and a unique family $\{F_i:i<\alpha\}$ of perfect dilators such that
\[F = \sum_{\xi<\alpha}F_\xi.\]
(A dilator is \emph{perfect} if it is additively indecomposable, i.e., whenever $F_i = F' + F''$, then either $F' = 0$ or $F'' = 0$.) Moreover, the perfect decomposition theorem asserts that if $G = \sum_{\xi < \beta} G_\xi$ is another sum of perfect dilators and $T:F\to G$ is a natural transformation, then there is a unique embedding $h : \alpha\to\beta$ and a unique family $\{T_\xi:\xi<\alpha\}$ of natural transformations from $F_\xi$ to $G_{h(\xi)}$ such that $T = \sum_{\xi < \alpha} T_\xi$ (where this sum is defined in the natural way).


Now we construct an embedding from $D^*$ to $\hat D$. We apply the perfect decomposition theorem to both $D^*$ and $\hat D$ and write:
\[D^* = \sum_{\zeta<\alpha}D^*_\zeta\;\;\;\;\;\;\hat D = \sum_{\delta<\beta}\hat D_\delta.\]

Observe that for any $\gamma<\alpha$ the dilator  $S_{\gamma}=\sum_{\zeta<\alpha'}D^*_\zeta$ is a subdilator of some $\sum_{i<n}D_i$. Since each $\sum_{i<n}D_i$ is a $T$-provable recursive dilator, for each $\gamma<\alpha$ the dilator  $S_{\gamma}$ is embeddable into $\hat D$. Combining this with the perfect decomposition theorem we see that for any $\gamma$ there is a function $f_\gamma\colon \gamma\to \beta$ such that for any $\zeta<\gamma$ the perfect dilator $D^*_\zeta$ is embeddable into $\hat D_{f_\gamma(\zeta)}$.  

We define the function $f\colon \alpha\to \beta$:
\[f(\zeta)=\min \{f_\gamma(\zeta)\mid \alpha>\gamma>\zeta\}.\]
Since all $f_\gamma$ were strictly increasing, $f$ is strictly increasing. And since for any $\gamma<\alpha$ and $\zeta<\gamma$ the dilator $D^*_\zeta$ was embeddable into $\hat D_{f_\gamma(\zeta)}$, we see that for any $\zeta<\alpha$ the dilator $D^*_\zeta$ is embeddable into $\hat D_{f(\zeta)}$.

For each $\zeta<\alpha$ we fix an embeeding $e_\zeta\colon D^*_\zeta\to \hat D_{f(\zeta)}$. We define the embedding $e\colon D^*\to \hat D$ to be the sum $\sum_{\zeta<\alpha}e_\zeta$. This proves the existence of $|T|_{\Pi^1_2}$.

For the ``moreover'' part of the theorem, we have to modify the definition of $D^*$ slightly.
We first observe that, with the definition given, $D^*$ is recursively enumerable if $T$ is. The remainder of the proof consists in a variant of the usual technique for replacing recursively enumerable structures by a recursive isomorphic copy.
This is done as follows: regarding $D^*$ as a denotation system for ordinals, let $t$ be a term for an ordinal in $D^*$. Thus, $t$ is a term in $D_i$, for some $i \in\mathbb{N}$. We define a dilator $D^{**}$ consisting of pairs $(t,p)$, where $t$ is as above, and $p$ is the least $T$-proof witnessing the fact that $D_i$ is indeed the $i$th provable recursive dilator of $T$. All comparison rules for terms in $D^{**}$ are the same as for $D^{*}$. It should be clear that $D^{**}$ is as desired.
\endproof

\begin{corollary}\label{CorollaryPi12Dilator}
Suppose $T$ is $\Pi^1_2$-sound. Then $|T|_{\Pi^1_2} = \sum_{D\in\mathcal{D}} D$, where $\mathcal{D}$ is an enumeration of all the provable recursive dilators of $T$ in order type $\omega$. Moreover, up to bi-embeddability, $|T|_{\Pi^1_2}$ does not depend on the enumeration used.
\end{corollary}
\proof
This was part of the proof of Theorem \ref{TheoremPi12Dilator}. The fact that $|T|_{\Pi^1_2}$ does not depend on the enumeration used follows from the fact that $\mathcal{D}$ is closed under finite sums.
\endproof

\begin{corollary}
Suppose that $T$ is $\Pi^1_2$ sound and recursively enumerable. Then, for every provable dilator $D$ of $T$, there is a recursive natural transformation from $D$ to $|T|_{\Pi^1_2}$.
\end{corollary}
\proof
This was obtained during the ``moreover'' part of the proof of Theorem \ref{TheoremPi12Dilator}.
\endproof

\begin{theorem}
  For any theory $T$ extending $\aca$ the statement ``$|T|_{\Pi^1_2}$ is a dilator.'' is equivalent to the scheme $\Pi^1_2\textsf{-RFN}(T)$ of uniform $\Pi^1_2$-reflection for $T$:
 \begin{equation} \label{RFN_inst}\forall x(\mathsf{Prv}_T(\ulcorner \varphi(\dot x)\urcorner)\to \varphi(x)),\text{ for }\varphi\in \Pi^1_2.\end{equation}
\end{theorem}
\proof
We reason in $\aca$.

First we assume $\Pi^1_2\textsf{-RFN}(T)$ and prove that $|T|_{\Pi^1_2}$ is a dilator. Indeed the latter is equivalent to the assertion that for any recursive $D$ if $T\vdash \text{``$D$ is a dilator''}$ then $D$ is a dilator. This is implied by $\Pi^1_2\textsf{-RFN}(T)$ since formulas ``$D$ is a dilator'' are $\Pi^1_2$.

Now we assume that $|T|_{\Pi^1_2}$  is a dilator and prove an instance of (\ref{RFN_inst}). Since $\varphi(x)$ is $\Pi^1_2$ and $\mathrm{DIL}$ is $\Pi^1_2$-complete in $\aca$, there are recursive dilators $D_x$ such that $\varphi(x)\mathrel{\leftrightarrow} \text{``$D_x$ is a dilator''}$ and this equivalence is $\aca$-provable. To finish the proof we assume that $\varphi(x)$ is $T$-provable and claim that $\varphi(x)$ is true. $T$ proves that $D_x$ is a dilator. Thus, being embeddable into $|T|_{\Pi^1_2}$, $D_x$ is a dilator. Therefore $\varphi(x)$ is true.\endproof

In \cite{PW21}, the authors defined the notion of a \emph{proof-theoretic dilator}, though the definition was only extensional (i.e., only defined for wellorders) and not functorial. It is natural to wonder how this definition relates to the present notion of the $\Pi^1_2$-rank of a theory. We will show that $|T|_{\Pi^1_2}$ is extensionally equal to the proof-theoretic dilator of $T$. First, we need a lemma.
\begin{lemma}\label{LemmaImplicationDilator}
Suppose $a$ and $b$ are countable linear orders. Then, we can effectively find a pre-dilator $D_{a\to b}$ such that
\begin{enumerate}
\item there is an embedding $e: b \to D_{a\to b}(a)$;
\item $D_{a\to b}$ is a dilator if and only if the implication $\wo(a)\to\wo(b)$ holds;
\item if $D_{a\to b}$ is not a dilator, $D_{a\to b}(x)$ is illfounded  if and only if $a$ embeds into $x$.
\end{enumerate}
\end{lemma}
\proof
We define $D_{a\to b}$. $D_{a\to b}(x)$ is the Kleene-Brouwer ordering of the tree $T_{a,b}(x)$ of attempts to simultaneously build:
\begin{enumerate}
\item an infinite descending chain through $b$, and
\item an embedding from $a$ into $x$.
\end{enumerate}

Let us now define $T_{a,b}(x)$ more formally. We suppose $a$ and $b$ are orders whose domain is a subset of $\omega$. The nodes of $T_{a,b}(x)$ are triples $\langle n, f,g\rangle$, where $n$ is a natural number, $f\colon n\to b$ is a descending sequence in $b$ and $g\colon a\upharpoonright n \to x$ is an order-preserving function. The root of the tree is $\langle 0,\emptyset,\emptyset\rangle$. The children of a node $\langle n, f,g\rangle$ are nodes $\langle n+1, f',g'\rangle$ such that $f'\upharpoonright n=f$ and $g'\upharpoonright n=g$. We define the comparison of two children $\langle n+1, f_1,g_1\rangle\ne \langle n+1, f_2,g_2\rangle$ of the same node $\langle n, f,g\rangle$: If $f_1(n)\ne f_2(n)$ then we compare the children according to the comparison of $f_1(n)$ and $f_2(n)$ as natural numbers; if $f_1(n)=f_2(n)$ then we compare the children according to the comparison of $g_1(n)$ and $g_2(n)$ as elements of $x$ (note that we could have $f_1(n)=f_2(n)$ only if $n\in \dom(a)$). It is easy to see that when $x$ is a well-order the order on the children of any node in $T_{a,b}(x)$ is a well-order. 

From the construction it is easy to see that $D_{a\to b}$ is is a pre-dilator.

Let us show that $D_{a\to b}$ is a dilator if and only if the implication $\wo(a)\to\wo(b)$ holds. If $b$ is wellfounded, then $T_{a,b}(x)$ is always wellfounded, and thus $D_{a\to b}(x)$ is wellfounded for all well-orders $x$. If $b$ is illfounded, then for well-orders $x$ the tree $T_{a,b}(x)$ (and thus the order $D_{a\to b}(x)$) is illfounded if and only if there is an embedding from $a$ into $x$. 

It remains to prove that there is an embedding from $b$ into $D_{a\to b}(x)$ if $x$ embeds into $a$. In this case, it suffices to see that (i) $b$ embeds into the Kleene-Brouwer order on the tree of all descending chains through $b$ (that is in $D_{\emptyset\to b}(\emptyset)$) (ii) this order embeds into $D_{a,b}(x)$. It is trivial to see that (ii) holds: if we have an embedding $u\colon a\to x$, then we embed $D_{\emptyset\to b}(\emptyset)$ to $D_{a,b}(x)$ by mapping $\langle n,f,\emptyset\rangle$ to $\langle n,f,u\upharpoonright n\rangle$.

Let us define an embedding $e\colon b \to D_{\emptyset\to b}$ and prove (i). We enumerate all elements of $b$ in order-type $\omega$, say by $b_0,b_1,\ldots$ (if $b$ is finite, this sequence is finite). We define values $e(b_i)$ by induction on $i$.  If $b_i=\max_b\{b_j \mid j\le i\}$, then $e(b_i)=\langle 1,f,\emptyset\rangle$, where $f(0)=b_i$. Otherwise we consider $b_l=\min_b\{b_j>_bb_i \mid j< i\}$ and $e(b_l)=\langle n,f,\emptyset\rangle$, we put $e(b_i)=\langle n+1,f',\emptyset\rangle$, where $f'$ extends $f$ by $f'(n)=b_i$. A straightforward verification shows that $e$ indeed is an embedding.
\endproof

\begin{theorem} \label{TheoremPi12DilatorandP11Ordinal}
Suppose that $T$  is $\Pi^1_2$-sound and extends $\aca$. and let $D = |T|_{\Pi^1_2}$. Then, for every recursive wellorder $\alpha$, we have
\[D(\alpha) = |T + \wo(\alpha)|_{\Pi^1_1}.\]
\end{theorem}
\proof
We first show that $D(\alpha) \leq |T+\wo(\alpha)|_{\Pi^1_1}$. By Corollary \ref{CorollaryPi12Dilator}, 
\[D(\alpha) = \sum_{i\in\mathbb{N}} D_i(\alpha),\]
where $D_0, D_1, \hdots$ enumerates all the provable dilators of $T$. For any fixed $j$, $D_0 + D_1 + \dots + D_j$ is a provable dilator of $T$, so
\[T + \wo(\alpha) \vdash \wo(D_0(\alpha) + \dots + D_j(\alpha)).\]
It follows that 
\[D_0(\alpha) + \dots + D_j(\alpha)  < |T+\wo(\alpha)|_{\Pi^1_1}\]
for each $j\in\mathbb{N}$ and thus
\[\sum_{i\in\mathbb{N}}D_i(\alpha) \leq |T + \wo(\alpha)|_{\Pi^1_1}.\]
We now prove that $|T+\wo(\alpha)|_{\Pi^1_1} \leq D(\alpha)$. Suppose that
\[T + \wo(\alpha) \vdash \wo(\beta).\]
Then, we have
\[T \vdash \wo(\alpha)\to\wo(\beta).\]
By Lemma \ref{LemmaImplicationDilator}, 
\[T\vdash \text{``$D_{\alpha\to \beta}$ is a dilator.''}\]
Moreover by Lemma \ref{LemmaImplicationDilator}, there is an embedding $e: \beta\to D_{\alpha\to \beta}(\alpha)$. Since $D_{\alpha\to\beta}$ is a provable dilator of $T$, there is a  natural transformation from $D_{\alpha\to\beta}$ to $D$, so we get an embedding 
$e': \beta\to D(\alpha)$, thus proving the inequality.
\endproof

\section{The $\Pi^1_2$-soundness ordinal}\label{SectUnsoundness}
The purpose of this section is to consider $\Pi^1_2$ analyses beyond $\Pi^1_2$-sound theories $T$. For such a $T$, we cannot expect to obtain a dilator bounding all the $T$-provable dilators. However, we can distil from the proof of existence of $\Pi^1_2$-norms (\mbox{Theorem \ref{TheoremPi12Dilator}}) an  alternative measure which will provide useful information about $T$. Since the $\Pi^1_2$ consequences of $T$ cannot be relied upon, we are forced to retreat back into the realm of ordinals.

\begin{definition}\label{Defo12}
Let $T$ be a theory. We define the \emph{$\Pi^1_2$-soundness ordinal} of $T$ by
\[o^1_2(T) = \min\Big\{\alpha: \Big(\sum_{D\in \mathcal{D}}D\Big)(\alpha)\text{ is illfounded}\Big\},\]
where $\mathcal{D}$ is an enumeration of all the $T$-provable recursive dilators in order type $\omega$.
\end{definition}
According to Corollary \ref{CorollaryWellDefined} below, the definition of $o^1_2(T)$ does not depend on the enumeration of $\mathcal{D}$.
If $o^1_2(T)$ is undefined, we may write $o^1_2(T) = \infty$. If $o^1_2(T)  < \infty$, we may write $o^1_2(T) \in\Ord$.

\begin{lemma}
Suppose $T$ is $\Pi^1_2$-unsound. then $o^1_2(T)$ is the least ordinal $\alpha$ such that for some provable recursive dilator $D$ of $T$, $D(\alpha)$ is illfounded.
\end{lemma}
\proof
Letting $\mathcal{D}$ be an enumeration of all the provable dilators of $T$ in order type $\omega$, the lemma
follows from the observation that for each $\alpha$,
\[\Big(\sum_{D\in \mathcal{D}}D\Big)(\alpha) = \sum_{D\in\mathcal{D}}(D(\alpha))\]
is a sum of linear orders in order-type $\omega$, so it is illfounded if and only if one of the summands is illfounded.
\endproof

\begin{corollary}\label{CorollaryWellDefined}
$o^1_2(T)$ does not depend on the enumeration of $\mathcal{D}$ chosen.
\end{corollary}
\proof
Immediate.
\endproof

The ordinal $o^1_2(T)$ is a measure of how close $T$ is to being $\Pi^1_2$-sound. 
The main thesis of this section is that $o^1_2(T)$ can be used to obtain useful information about $T$. We will classify theories $T$  into four categories based on the value $o^1_2(T)$. The classification is very natural, and
we will see that for recursively enumerable extensions of $\aca$, theories in each category share similar  properties. 

\begin{definition}
Let $T$ be a theory. We say that:
\begin{enumerate}
\item $T$ is in Category $A$ if $o^1_2(T) = 0$;
\item $T$ is in Category $B$ if $0 < o^1_2(T)< \omega_1^{ck}$;
\item $T$ is in Category $C$ if $\omega_1^{ck}\leq  o^1_2(T) < \infty$;
\item $T$ is in Category $D$ if $o^1_2(T) = \infty$.
\end{enumerate}
\end{definition}
We will see below that only some values in Categories B and C can be attained.
An interesting problem is to
characterize the ordinals $\alpha$ which are of the form $o^1_2(T)$ for some recursively enumerable extension of $\aca$. This is the \emph{Spectrum Problem} for $\Pi^1_2$-soundness. In the following sections, we will come back to this problem and solve it for theories in Category B, as well as for theories in Category C whose $\Pi^1_2$-soundness ordinal is admissible.

First, we prove an optimal upper bound for the ordinals $o^1_2(T)$:
\begin{theorem}\label{TheoremDelta12}
Let $\delta^1_2 = \sup\{\alpha: \alpha$ is the length of a $\Delta^1_2$ wellordering of $\mathbb{N}\}$. Then
\begin{align*}
\delta^1_2 &= \text{$\sup\{ o^1_2(T) \in \Ord:$ $T$ is a recursively enumerable extension of $\eca\}$.}\\
&= \text{$\sup\{ o^1_2(T)\in\Ord:$ $T$ is a $\Sigma^1_2$ extension of $\eca\}$.}
\end{align*}
\end{theorem}
\proof 
We will use the fact that $\delta^1_2$  is the least ordinal $\sigma$ such that
\[L_\sigma\prec_1 L.\]
We refer the reader to Barwise \cite{Ba75} for a proof.

First we show that $\delta^1_2$ is an upper bound. Let $T$ be a $\Sigma^1_2$ extension of $\eca$ and suppose that
\[D(\alpha) \text{ is illfounded}\]
for some $\alpha$ and some provable dilator $D$ of $T$.
Since $T$ is $\Sigma^1_2$, there is a $\Sigma_1$ formula $\psi$ in the language of set theory such that $\varphi$ is an axiom of $T$ precisely when
\[L_{\delta^1_2}\models\psi(\ulcorner \varphi \urcorner)\]
holds. 
Given an ordinal $\gamma$, let $T^\gamma$ denote the theory whose axioms are all the formulas $\varphi$ such that 
\[L_\gamma\models \psi(\ulcorner \varphi \urcorner).\]
Since $\psi$ is $\Sigma_1$, it follows that $T^\gamma$ is a subtheory of $T$ whenever 
 $\gamma<\delta^1_2$. 
 For $\gamma\geq\delta^1_2$, we have $T^\gamma  = T$.

Let $\theta$ be the $\Sigma_1$ sentence in language of set theory expressing the existence of some ordinal $\gamma$, some recursive pre-dilator $D$, and some ordinal $\alpha$ such that:
\begin{enumerate}
\item $T^\gamma$ proves that $D$ is a dilator; and
\item there is an infinite descending chain through $D(\alpha)$.
\end{enumerate}
This is a true $\Sigma_1$ assertion in the language of set theory. 
By Shoenfield absoluteness it holds in $L$. Since $L_{\delta^1_2} \prec_1 L$, we have 
\[L_{\delta^1_2}\models\theta.\] 
Hence, we get some $\gamma<\delta^1_2$, a pre-dilator $D$ and some ordinal $\alpha<\delta^1_2$ such that $T^\gamma$ proves that $D$ is a dilator, but $D(\alpha)$ has an infinite descending chain in $L_{\delta^1_2}$. Since $T^\gamma$ is a subtheory of $T$, $T$ proves that $D$ is a dilator.
This shows that $\delta^1_2$ is an upper bound.\\

For the lower bound, let $\delta < \delta^1_2$ be arbitrary. We find a recursively enumerable extension of $T$ such that $\delta < o^1_2(T)$.  Since $\delta <\delta^1_2$, it follows that 
\[L_\delta\not\prec_1 L,\]
so there is $\bar\delta$ with $\delta<\bar\delta$ and a $\Sigma_1$ formula $\varphi$ such that $L_{\delta}\not\models \varphi$, but $L_{\bar\delta}\models\varphi$.
Without loss of generality, we may assume that $L_{\bar\delta}$ is a limit of $\gamma$ such that $L_\gamma\models\Pi^1_2$-CA$_0$ (this is because we can replace $\varphi$ with the conjunction of $\varphi$ and ``there are arbitrarily $\gamma$ such that $L_\gamma\models\Pi^1_2$-CA$_0$''). 
Let $T$ be the theory consisting of
\begin{enumerate}[resume]
\item all axioms of $\aca$;
\item an axiom asserting that there is no countably-coded well-founded model of $\varphi$ + V = L.
\end{enumerate}
We claim that $\delta < o^1_2(T)$. To see this, choose some $\gamma$ such that $\delta < \gamma < \bar\delta$ and 
\[L_\gamma\models\Pi^1_2\text{-CA}_0.\]
Then, $T$ is $\Pi^1_2$-sound in $L_\gamma$. Since the proof of Theorem \ref{TheoremPi12Dilator} goes through in $\Pi^1_2$-CA$_0$ (and much weaker theories), we have
\[L_\gamma\models o^1_2(T) = \infty.\]
Since $\gamma$ is recursively inaccessible, $L_\gamma$  is correct about wellfoundedness, so it follows that  $\gamma < o^1_2(T)$, as desired.
\endproof
In a later section, we will refine the ideas behind (the second part of) the proof of Theorem \ref{TheoremDelta12} in order to produce theories with a specific admissible $\Pi^1_2$-soundness ordinal.
\begin{remark}
In the statement of Theorem \ref{TheoremDelta12}, we could replace $\eca$ by any $\Pi^1_2$-sound theory which has a transitive model. In particular, we could replace it by $\aca$ or by $\zfc$ (under suitable set-theoretic assumptions).
\end{remark}

We can now derive a characterization of the theories in Category D.
\begin{proposition}
The following are equivalent:
\begin{enumerate}
\item $T$ is in Category D;
\item $\delta^1_2 \leq o^1_2(T)$;
\item $T$ is $\Pi^1_2$-sound.
\end{enumerate}
\end{proposition}
\proof
This is immediate from  Corollary \ref{CorollaryPi12Dilator} and Theorem \ref{TheoremDelta12}.
\endproof

We state two theorems characterizing soundness properties for theories in terms of their $\Pi^1_2$-soundness ordinals. The following theorem asserts that the $\Pi^1_1$-sound, recursively enumerable extensions of $\aca$ are precisely those in Categories B, C, and D. Moreover, the condition $o^1_2(T) < |T|_{\Pi^1_1}$ is enough to guarantee that $T$ belongs to Category A.
\begin{theorem}\label{TheoremPi12Zero}
Suppose $T$ is a recursively enumerable extension of $\aca$. The following are equivalent:
\begin{enumerate}
\item $T$ is $\Pi^1_1$-sound;
\item $0 < o^1_2(T)$;
\item $|T|_{\Pi^1_1} \leq o^1_2(T)$.
\end{enumerate}
\end{theorem}
\proof
If $T$ is not $\Pi^1_1$-sound, then there is some recursive illfounded linear ordering $a$ such that 
\[T\vdash \wo(a).\]
Letting $C_a$ be the constant pre-dilator with value $a$, we have 
\[T\vdash \text{``$C_a$ is a dilator.''}\]
In particular, $T$ proves that $C_a(0)$ is wellordered, while in reality it is not, so $o^1_2(T) = 0$. 
Suppose now that $T$ is $\Pi^1_1$-sound and that $\alpha<|T|_{\Pi^1_1}$, so that
\[T\vdash\wo(\alpha).\]
Suppose moreover that 
\[T\vdash \text{``$D$ is a dilator''}\]
for some recursive pre-dilator $D$. Then,
\[T\vdash \wo(D(\alpha)).\]
Since $T$ is $\Pi^1_1$-sound, then $D(\alpha)$ really is wellordered. Hence, $\alpha < o^1_2(T)$. We conclude $|T|_{\Pi^1_1} \leq o^1_2(T)$, as desired.
\endproof

Let $\mathsf{Bool}(\Pi^1_1)$ denote the collection of Boolean combinations of $\Pi^1_1$ sentences.  The following theorem asserts that the $\mathsf{Bool}(\Pi^1_1)$-sound extensions of $\aca$ are precisely those in Categories C and D.
\begin{theorem}\label{TheoremPi12NonRecursive}
Suppose $T$ is a recursively enumerable extension of $\aca$. The following are equivalent:
\begin{enumerate}
\item \label{TheoremPi12NonRecursive1} $o^1_2(T)$ is non-recursive;
\item 
\label{TheoremPi12NonRecursive3}
$T$ is $\mathsf{Bool}(\Pi^1_1)$-sound.
\end{enumerate}
\end{theorem}
\proof
We prove that \eqref{TheoremPi12NonRecursive3} implies \eqref{TheoremPi12NonRecursive1}.
Suppose that $o^1_2(T)$ is recursive. Let $a$ be a recursive wellorder such that for some recursive pre-dilator $D$, the following hold:
\begin{enumerate}
\item $T\vdash $``$D$ is a dilator,''
\item $D(a)$ is illfounded.
\end{enumerate}
Then, 
\[T\vdash \wo(a)\to\wo(D(a)),\]
which is a false implication of $\Pi^1_1$ sentences. 

We now prove that \eqref{TheoremPi12NonRecursive1} implies \eqref{TheoremPi12NonRecursive3}. We first suppose that $T$ is not $\Pi^1_1\to\Pi^1_1$-sound, so there are $\Pi^1_1$ sentences $\phi,\psi$ such that $\phi\to\psi$ is false, but 
\[T\vdash\phi\to\psi.\]
By the completeness of $\wo$ (which is provable in $\aca$), we may assume that $\phi$ is of the form $\wo(a)$ and $\psi$ is of the form $\wo(b)$, for some recursive wellorder $a$ and some recursive illfounded linear order $b$.
Thus,
\[T\vdash \wo(a)\to\wo(b).\]
By Lemma \ref{LemmaImplicationDilator},
\[T\vdash \text{``$D_{a\to b}$ is a dilator.''}\]
By Lemma \ref{LemmaImplicationDilator}, there is an embedding from $b$ to $D_{a\to b}(a)$. By assumption, $b$ is illfounded. Hence $D_{a\to b}(a)$ is illfounded. By assumption, $a$ is a wellorder, and it is recursive, so we have $o^1_2(T)\leq a < \omega_1^{ck}$, as desired.

We have shown that $o^1_2(T)$ is non-recursive if and only if $T$ is sound for sentences of the form $\Pi^1_1\to\Pi^1_1$. However, every Boolean combination of $\Pi^1_1$ sentences can be written in the form
\begin{equation}\label{eqBooleanPi11}
\phi_0 \wedge \phi_1 \wedge \dots \wedge \phi_n
\end{equation}
where each $\phi_i$ is of the form  $\Pi^1_1\to\Pi^1_1$. This can be shown by considering the class of all formulas 
logically equivalent to a formula of the form \eqref{eqBooleanPi11} and observing that this class is closed under conjunctions and complements.
Hence, if $T$ is sound for implications between $\Pi^1_1$ sentences, then it is sound for Boolean combinations of $\Pi^1_1$-sentences. 
\endproof
By putting together the last three results, one obtains Theorem \ref{thmxB} from the introduction.

\section{Theories with recursive $\Pi^1_2$-soundness ordinal}\label{SectRecursive}
In this section, we study theories in Category B and their $\Pi^1_2$-soundness ordinals. We shall obtain a solution to the  Spectrum Problem for these theories. First, a lemma:
\begin{lemma}\label{LemmaEpsilon}
Suppose $T$ is a $\Pi^1_1$-sound, recursively enumerable extension of $\aca$. Then, $o^1_2(T)$ is of the form $\varepsilon_\alpha$ for some $\alpha$.
\end{lemma}
\proof
Suppose towards a contradiction that for some $\alpha$, we have
\[\varepsilon_\alpha < o^1_2(T) < \varepsilon_{\alpha+1}.\]
This means that there is a pre-dilator $D$ such that $T$ proves that $D$ is a dilator and a least ordinal $\gamma$ with  such that  $D(\gamma)$ is illfounded and moreover $\gamma$ satisfies $\varepsilon_\alpha < \gamma< \varepsilon_{\alpha+1}$. Since $\varepsilon_\alpha < \gamma< \varepsilon_{\alpha+1}$, we have
\[\varepsilon_\alpha < \gamma < \omega^{\omega^{\iddots^{\varepsilon_{\alpha}+1}}}\]
for some natural number $n$. Let $F$ be a dilator such that
\[F(x) = \omega^{\omega^{\iddots^{x+1}}}\]
for all $x$ and such that $F$ is a dilator provably in $\aca$. Such an $F$ exists by a theorem of Girard whereby $\aca$ is equivalent over $\rca$ to the statement that
\[x\mapsto \omega^x\]
preserves wellfoundedness. But then we have
\[T\vdash \text{``$D\circ F$ is a dilator''}\]
and $D\circ F(\varepsilon_\alpha)$ is illfounded, contradicting the choice of $\gamma$.
\endproof

\begin{lemma}\label{LemmaDisjOrdinal}
Suppose $a$ and $b$ are countable linear orders. 
Then, we can effectively and uniformly find a linear order $l_{a\vee b}$ such that, provably in $\aca$,
\begin{enumerate}
\item if either of $a$ or $b$ is wellfounded, then $l_{a\vee b}$ is wellfounded;
\item if $a$ is wellfounded and $b$ is illfounded, then $a$ embeds into $l_{a\vee b}$;
\item if $a$ is illfounded and $b$ is wellfounded, then $b$ embeds into $l_{a\vee b}$.
\end{enumerate}
Moreover, for each countable linear order $b$, there is a pre-dilator $F$ such that 
\[F(a) = l_{a\vee b}\]
for all $a$.
\end{lemma}
\proof
The first part of the lemma is stated and proved in \cite{PW21}. The ``moreover'' part follows from the uniformity of the construction.
\endproof

The following theorem is the solution to the Spectrum Problem for theories in Category B.
\begin{theorem}\label{TheoremSpectrumRecursive}
The recursive ordinals of the form $o^1_2(T)$ for some recursively enumerable extension of $\aca$ are precisely the ordinals of the form $\varepsilon_\alpha$, for $\alpha<\omega_1^{ck}$.
\end{theorem}
\proof
Fix $\alpha<\omega_1^{ck}$ and an ordinal notation system for $\varepsilon_\alpha$.
Let $T$ be the theory consisting of the following axioms:
\begin{enumerate}
\item \label{eqSpectrumRecursiveT1} $\aca$;
\item \label{eqSpectrumRecursiveT3} ``$\varepsilon_\alpha$ is illfounded.''
\end{enumerate}
Let $\omega^\star$ be an infinite descending chain and let $D_{\varepsilon_\alpha\to\omega^\star}$ be as in Lemma \ref{LemmaImplicationDilator}. Notice that, by Lemma \ref{LemmaImplicationDilator}, the theory $T$ proves that $D_{\varepsilon_\alpha\to\omega^\star}$ is a dilator and the order $\omega^\star$ embeds into $D_{\varepsilon_\alpha\to\omega^\star}(\varepsilon_\alpha)$. Hence $o^1_2(T)\le \varepsilon_\alpha$.

To prove that $o^1_2(T)=\varepsilon_\alpha$, we consider an arbitrary $\beta<\varepsilon_\alpha$ and a $T$-provable recursive dilator $D$ and prove that $D(\beta)$ is well-ordered. It is a folklore result that $\aca+\wo(\beta)$ has as its proof-theoretic ordinal $\varepsilon_\gamma$ for the least $\gamma$ such that $\varepsilon_\gamma>\beta$ (this can be shown by a slight modification of the ordinal analysis of $\aca$). Clearly $\varepsilon_\gamma\le \varepsilon_\alpha$. Observe that 
\[\aca+\wo(\beta)\vdash \lnot \wo(\varepsilon_\alpha)\to \wo(D(\beta)).\] 
Hence, in the notation of Lemma \ref{LemmaDisjOrdinal}, we have $\aca+\wo(\beta)\vdash\wo(l_{\varepsilon_\alpha\lor D(\beta)})$. Therefore $\mathsf{ot}(l_{\varepsilon_\alpha\lor D(\beta)})<\varepsilon_\alpha$, which is only possible if $D(\beta)$ is well-ordered, since otherwise by Lemma \ref{LemmaDisjOrdinal}, $\varepsilon_\alpha$ would be embeddable into $l_{\varepsilon_\alpha\lor D(\beta)}$.
\endproof

The proof of Theorem \ref{TheoremSpectrumRecursive}  illustrates what a theory $T$ must look like in order for $o^1_2(T)$ to take a specific recursive value. Thus, just like the value $o^1_2(T)$ is a measure of how close $T$ is to being $\Pi^1_2$-sound,  \emph{recursive} values of $o^1_2(T)$ are a measure of how close $T$ is to being $\mathsf{Bool}(\Pi^1_1)$-sound.\\

\noindent We now turn to theories in Category C.

\section{Theories with admissible $\Pi^1_2$-soundness ordinal}\label{SectAdmissible}
In this section, we study theories in Category C and their $\Pi^1_2$-soundness ordinals. We will obtain a solution to the Spectrum Problem for theories with admissible $\Pi^1_2$-soundness ordinal.

Below, for an ordinal $\alpha$, we denote by $\alpha^+$ the smallest admissible ordinal larger than $\alpha$.

Recall that a $\beta$-model $M$ is an $\omega$-model of second-order arithmetic that satisfies all true $\Sigma^1_1$-sentences with parameters from $M$. An important feature of $\beta$-models is that they are correct about the well-foundedness of sets. For a $\beta$-model $M$ we denote by $o(M)$ the least ordinal such that for any  well-ordering $a$ inside $M$ we have $a<o(M)$.

Recall that for a set of naturals $X$ we denote by $\omega_1^X$ the first admissible ordinal relative to $X$ (alternatively $\omega_1^X$ could be defined as the supremum of order types of $X$-recursive well-orderings). Notice that for a $\beta$-model $M$ and a set $X \in M$ we always have $\omega_1^X\le o(M)$.

\begin{lemma} \label{LemmaSoundOmegaModel}
For any admissible ordinal $\gamma$ there is a $\beta$-model $M\models\aca$ and a set $G\in M$ such that $o(M)=\omega_1^G=\gamma$.
\end{lemma}
\proof
Let $\mathbb{S}$ be the partial order of Sacks \cite{Sa76} consisting of all hyperarithmetically pointed perfect trees $t$ which belong to $L_\gamma$. Conditions are subtrees $t\in L_\gamma$ of $2^{<\omega}$ with the property that $t$ is hyperarithmetical in every path through $t$.
Let $G\subset\mathbb{S}$ be sufficiently generic (it suffices that $G$ have nonempty intersection with all subsets of $\mathbb{S}$ definable over $L_{\gamma}$). By a theorem of Sacks \cite{Sa76} we have
$L_\gamma[G] \models\kp$
and 
\[\omega_1^G = \gamma.\]

To construct $M$ we apply Corollary VII.2.12 of Simpson \cite[Corollary~VII.2.12]{Simpson}, from which it follows that there exists a $\beta$-model $M$ such that $G\in M$ and for all $X\in M$ we have $\mathcal{O}^X\le_T\mathcal{O}^G$. 

It follows that $\omega_1^X\leq\omega_1^G$ for all $X\in M$. Indeed, suppose otherwise that $\omega_1^G < \omega_1^X$ for some $X\in M$. Without loss of generality (e.g., by replacing $X$ with $(X,G)$ if necessary), we may assume that $G$ is recursive in $X$.
By a theorem of Spector \cite{Sp55} (see also Sacks \cite[II.7.6]{Sa90}), we have 
\[\omega_1^G < \omega_1^X \text{ and } G \leq_{hyp} X \text{ imply } \mathcal{O}^G \leq_{hyp} X,\]
so that $\mathcal{O}^G \leq_{hyp} X$ and thus $\mathcal{O}^X\leq_{hyp} X$, which is impossible.
\endproof

\begin{lemma}\label{LemmaSigma11KPU}
For every $\Sigma^1_1$ sentence $\varphi$ in the language of set theory one can effectively and uniformly find a first-order formula $\varphi^*(X)$ in the language consisting of a binary predicate $<$ and a ternary predicate $X$, such that for every countable ordinal $\alpha$, the following are equivalent:
\begin{enumerate}
\item $(L_\alpha,\in) \models \varphi$ and $\alpha$ is a limit ordinal; and
\item there is $X\subset \alpha\times \alpha\times \alpha$ such that $(\alpha,<, X)\models\varphi^*(X)$.
\end{enumerate} 
\end{lemma}
\proof
Recall that second-order logic allows quantification over $n$-ary relations.
Letting $X_i = \{(x,y): (i,x,y) \in X\}$,
the formula $\varphi^*(X)$ is a formalization of the assertion that $\alpha$ is a limit ordinal,  where $X_0$ is a model of $V = L$, $X_1$ is a bijection between the ordinals of $X_0$ and $\alpha$, and $X_2$ is a subset of $X_0$ which witnesses $\varphi$ over the set coded by $X_0$.
\endproof

We recall the following definition: we say that $\alpha$ \emph{reflects} a formula $\varphi$ if 
\[L_\alpha\models \varphi \text{ implies } \exists \bar\alpha<\alpha\, L_{\bar\alpha}\models\varphi.\]
We say that $L_\alpha$ is \emph{parameter-free} (or \emph{lightface}) \emph{$\Sigma^1_1$-reflecting} if it reflects every $\Sigma^1_1$ sentence without parameters. We remark that in most cases, a countable ordinal which is parameter-free $\Sigma^1_1$-reflecting will have a bijection with $\mathbb{N}$ definable in a $\Sigma^1_1$ way (the first ordinal not satisfying this will be much greater than the least $\beta$-model of analysis) and thus will be $\Sigma^1_1$-reflecting with parameters as well.

\begin{theorem}\label{TheoremSpectrumAdmissible}
Let $\alpha$ be an admissible ordinal or a limit of admissibles. The following are equivalent:
\begin{enumerate}
\item $\alpha = o^1_2(T)$ for some recursively enumerable extension $T$ of $\aca$; and
\item $\alpha$ is not parameter-free $\Sigma^1_1$-reflecting.
\end{enumerate}
\end{theorem}
\proof
If $\alpha = o^1_2(T)$ then there is a recursive pre-dilator $D$ such that the following hold:
\begin{enumerate}
\item \label{condAdm1} $T\vdash \text{``$D$ is a dilator''}$;
\item \label{condAdm2} $D(\gamma)$ is wellfounded for every $\gamma<\alpha$;
\item \label{condAdm3} $D(\alpha)$ is illfounded.
\end{enumerate}
Since $D$ is recursive, it belongs to $L_{\omega+1}$. Since $\alpha$ is either admissible or a limit of admissible, condition \eqref{condAdm2} implies that $D(\gamma)<\alpha$ whenever $\gamma<\alpha$, so it can be expressed as a first-order assertion about $L_{\alpha}$. From this and the fact that $T$ is recursively enumerable, it follows that the conjunction of \eqref{condAdm1}--\eqref{condAdm3} can be expressed as a $\Sigma^1_1$ formula which holds of $L_\alpha$ but not of any $L_\gamma$ with $\gamma<\alpha$.\\

Conversely, suppose that $\alpha$ is not parameter-free $\Sigma^1_1$-reflecting and let $\varphi$ be a $\Sigma^1_1$ sentence such that $L_\alpha\models\varphi$ but $L_\gamma\not\models\varphi$ for all $\gamma<\alpha$. 
Let $\varphi^*(X)$ be the formula given by Lemma \ref{LemmaSigma11KPU} applied to $\varphi$. This is a first-order formula in which the relation symbol $X$ appears. 

Let $P = \{P(\xi):\xi\in\Ord\}$ be the $\beta$-pre-proof of the formula $\lnot\varphi^*(X)$ obtained from
the completeness theorem for $\beta$-logic (Theorem \ref{TheoremBetaCompleteness}), so that the following are equivalent:
\begin{enumerate}[resume]
\item $P(\xi)$ is wellfounded for every $\xi$; and
\item $\lnot\varphi^*(X)$ is $\beta$-valid.
\end{enumerate}
Since $\lnot\varphi^*(X)$ is not $\beta$-valid, $P(\xi)$ is not wellfounded for all $\xi$ and indeed $P(\xi)$ is illfounded precisely when $\alpha\leq\xi$.

Let $T$ be the theory
\[\aca + \text{``$P$ is a $\beta$-proof.''}\]
We claim that $o^1_2(T) = \alpha$.
Working in $T$, let $D$ be the pre-dilator which maps each $\gamma$ to the Kleene-Brouwer ordering of the proof tree $P(\gamma)$. Then, 
\[T\vdash \text{``$D$ is a dilator.''}\]
By choice of $P$, $P(\alpha)$ is illfounded, so $D(\alpha)$ is illfounded. Hence, $o^1_2(T)\leq\alpha$.

We need to show that if $\bar\alpha<\alpha$ and $D$ is a provable dilator of $T$, then $D(\bar\alpha)$ is wellfounded. Fix such an $\bar\alpha$ and
let $\gamma$ be the least admissible greater than $\bar\alpha$; thus, $\gamma$ is a successor admissible ($\gamma = \alpha$ is possible).
Consider the $\beta$-model $M$ provided by Lemma \ref{LemmaSoundOmegaModel}. Notice that since $o(M)=\gamma\le \alpha$, the formula $\lnot\varphi^*(X)$ holds for any well-ordering in the sense of $M$. Hence $P$ is a $\beta$-proof in $M$ and thus $D$ is a dilator in $M$. Thus from the perspective of $M$, the order $D(\bar\alpha)$ is wellfounded and hence $D(\bar\alpha)$ is well-ordered by correctness.
\endproof

\section{$\Sigma^1_2$-sound theories}\label{SectSigma}
In this section, we study the dual notion of $o^1_2(T)$. Namely, an ordinal which quantifies the complexity of true $\Sigma^1_2$ theorems provable by a theory.
\begin{definition}
Let $T$ be a $\Sigma^1_2$-sound extension of $\aca$. We define
\[s^1_2(T) = \sup\Big\{\min\big\{\alpha:  D(\alpha) \text{ is illfounded}\big\} : \text{$T \vdash $``$D$ is not a dilator''}\Big\}.\]
\end{definition}
Although we are interested in the case where $T$ is recursively enumerable, some of the arguments require that $s^1_2(T)$ be defined for more complicated theories $T$.

\begin{proposition}\label{Propo12Lesss12}
Let $T$ be a $\Sigma^1_2$-sound extension of $\aca$. Then, $s^1_2(T) \leq o^1_2(T)$.
\end{proposition}
\proof
Suppose towards a contradiction that $T$ is a $\Sigma^1_2$-sound extension of $\aca$ and $o^1_2(T) < s^1_2(T)$, witnessed by recursive pre-dilators $D_o$ and $D_s$, i.e., the following hold:
\begin{enumerate}
\item $T \vdash$ ``$D_o$ is a dilator,''
\item $T \vdash$ ``$D_s$ is not a dilator,''
\item \label{eqPropo12Lesss12} $T \vdash$ $\exists \alpha\, \wo(\alpha)\wedge \wo(D_o(\alpha))\wedge \lnot\wo(D_s(\alpha))$.
\end{enumerate}
By assumption, $o^1_2(T) < s^1_2(T)$, we get that if $s^1_2(T) \leq\alpha$, then $D_o(\alpha)$ is illfounded. In other words, if $\alpha$ is such that  $T\vdash \lnot\wo(D_s(\alpha))$, then $D_o(\alpha)$ is illfounded.
Hence, \eqref{eqPropo12Lesss12} is a false $\Sigma^1_2$ sentence, contradicting the choice of $T$.
\endproof

\begin{lemma}\label{Lemmas12upperbound}
Suppose $T$ is a $\Sigma^1_2$-sound extension of $\aca$. Suppose that $L_\alpha\models T$ for some recursively inaccessible $\alpha$ and $L_\sigma\prec_1 L_\alpha$ for some $\sigma < \alpha$. Then $s^1_2(T)\leq\sigma$.
\end{lemma}
\proof
Let $T'$ be the set of all $\Sigma^1_2$ consequences of $T$. Then, $T'$ is $\Sigma^1_2$-sound and $s^1_2(T) = s^1_2(T')$. Since $L_\sigma\prec_1 L_\alpha$, we have $L_\sigma\models T'$. Let $D$ be such that 
\[\text{$T \vdash $``$D$ is not a dilator.''}\]
The hypothesis implies that $\sigma$ is a limit of admissibles, so $L_\sigma \models $``$D$ is not a dilator,'' so there is an ordinal $\gamma<\sigma$ such that $D(\gamma)$ is illfounded. Thus, $s^1_2(T')\leq\sigma$.
\endproof

The following result gives a characterization of the set of ordinals of the form $s^1_2(T)$, for $\Sigma^1_2$-sound theories $T$ which extend $\Pi^1_1$-CA$_0$. Although the statement might seem like a trivial equivalence at first, it is not. For instance, it implies that $s^1_2(T)$ is always a limit of admissibles for such theories. In contrast to this, it is not hard to modify the earlier constructions to find an example of a recursively enumerable extension $T$ of $\Pi^1_1$ with $o^1_2(T) = \psi(\varepsilon_{\Omega+1})$.

\begin{theorem} \label{TheoremCharacterizations12}
Let $\alpha$ be an ordinal.
The following are equivalent:
\begin{enumerate}
\item \label{TheoremCharacterizations121} $\alpha= s^1_2(T)$ for some $\Sigma^1_2$-sound, recursively enumerable extension $T$ of $\Pi^1_1{-}\ca$; 
\item \label{TheoremCharacterizations122}$\alpha$ is a limit of admissibles and $\alpha$ is the least ordinal such that $L_\alpha\models S$, for some recursively enumerable set of $\Sigma_1$ sentences $S$ in the language of set theory.
\end{enumerate}
\end{theorem}
\proof
Suppose $\alpha$ is as in \eqref{TheoremCharacterizations122} and $S$ is the corresponding theory. We describe $T$. To each $\psi \in S$ we associate a $\beta$-pre-proof $P^\psi = \{P^\psi(\gamma):\gamma\in\Ord\}$ of the statement ``there is no transitive model of $\psi$.'' More specifically, for each ordinal $\gamma$, $P^\psi(\gamma)$ is  an attempted $\gamma$-proof of an appropriate formalization of ``it is not the case that $V = L$ and $\psi$ holds'' 
obtained via a Sch\"utte-type construction as in the  $\beta$-completeness theorem (Theorem \ref{TheoremBetaCompleteness}), so that, by choice of $\psi$, $P^\psi$ is not a $\beta$-proof.
Thus, for each $\psi\in S$ there is some $\alpha_\psi$ such that:
\begin{enumerate}
\item $P^\psi(\alpha_\psi)$ is an illfounded proof tree,
\item for each $\bar\alpha<\alpha_\psi$, $P^\psi(\bar\alpha)$ is a wellfounded proof tree,
\item $\sup_{\psi}\alpha_\psi = \alpha$.
\end{enumerate}

\begin{claim}
$\alpha$ is the strict supremum of $\{\alpha_\psi:\psi \in S\}$.
\end{claim}
\proof
By assumption, $L_\alpha\models S$ and $\psi$ is $\Sigma_1$, so there is $\bar\alpha<\alpha$ such that $L_{\bar\alpha}\models\psi$. Thus, $P^\psi(\bar\alpha)$ is illfounded, so $\alpha_\psi\leq\bar\alpha<\alpha$.
\endproof

We let $T$ be the theory consisting of the following sentences (suitably formalized):
\begin{enumerate}[resume]
\item $\Pi^1_1${-}$\ca$,
\item ``$P^\psi$ is not a $\beta$-proof,'' for each $\psi \in S$.
\end{enumerate}
Thus, $T$ is recursively enumerable. For each $\psi \in S$, we define a dilator $D_\psi$ which maps an ordinal $\gamma$ to the Kleene-Brouwer ordering on $P^\psi(\gamma)$. Hence, $D_\psi(\gamma)$ is illfounded if and only if $\alpha_\psi\leq\gamma$.
Since $T$ proves that $P^\psi$ is not a $\beta$-proof for each $\psi \in S$, $T$ proves that $D_\psi$ is not a dilator for each $\psi \in S$. Hence, we have $\alpha\leq s^1_2(T)$.

We need to show that $s^1_2(T)\leq\alpha$. 
Let $M$ be the $\omega$-model consisting of all sets of natural numbers which belong to $L_\alpha$.
\begin{claim}
$M \models T$.
\end{claim}
\proof
Since $\alpha$ is a limit of admissibles, we have $M\models\Pi^1_1{-}\ca$. If $\psi \in S$, then $P^\psi(\alpha_\psi)$ is illfounded by choice of $\alpha_\psi$. As we have seen, $\alpha$ is the strict supremum of $\{\alpha_\psi:\psi\in S\}$.
Thus, an infinite branch through $P^\psi(\alpha_\psi)$ is definable over $\alpha_\psi^+<\alpha$ and belongs to $L_\alpha$. 
\endproof
It follows that whenever
\[\text{$T \vdash $``$D$ is not a dilator,''}\]
we have
\[\text{$M \models $``$D$ is not a dilator,''}\]
and so there is a witness in $M$. Such a witness is hence
isomorphic to an ordinal ${<}\alpha$. We conclude that $s^1_2(T)\leq\alpha$. We have shown that \eqref{TheoremCharacterizations122} implies \eqref{TheoremCharacterizations121}.

For the converse, suppose $\alpha = s^1_2(T)$ for some $\Sigma^1_2$-sound recursively enumerable extension $T$ of $\Pi^1_1{-}\ca$. 
Let $\hat T$ be the set of all  sentences $\psi$ of the form ``$D_\psi$ is not a dilator'' for some recursive pre-dilator $D_\psi$ such that $\psi$ is provable in $T$.
Let 
\begin{align*}
S = \{\text{``there is $\alpha_\psi$ such that $D_\psi(\alpha_\psi)$ is illfounded''}: \psi\in \hat T\}.
\end{align*}
Since $T$ is recursively enumerable, $S$ is recursively enumerable. Let $\eta$ be least such that $L_\eta\models S$.
We claim that
\begin{enumerate}[resume]
\item\label{condChars121} $\eta = \alpha$,
\item\label{condChars122} $\alpha$ is  a limit of admissibles.
\end{enumerate}
For \eqref{condChars121}, we first observe that if $\psi \in \hat T$, then there is a witness $\alpha_\psi$ for $D_\psi$ not being a dilator, with $\alpha_\psi \in L_\eta$, so $\alpha\leq\eta$.
The following claim implies both the converse of \eqref{condChars121}, as well as \eqref{condChars122}.
\begin{claim}
Suppose $\gamma<\eta$. Then $\gamma^+ < \alpha$.
\end{claim} 
\proof
Suppose $\gamma<\eta$, so there is $\psi \in \hat T$ such that 
\[L_\gamma\models \forall \xi\in\Ord\, \wo(D_\psi(\xi)).\]
Consider the following sentence:
\begin{align*}
&\text{``it is not the case that $\kp$ holds, V = L, and}\\
&\, \text{$D_\psi(\xi)$ is illfounded for some $\xi$.''}
\end{align*}
Use the $\beta$-completeness theorem to find a $\beta$-pre-proof $P_\psi = \{P_\psi(\zeta):\zeta\in\Ord\}$ of (a suitable formalization of) the displayed sentence. If $\xi$ is such that $D_\psi(\xi)$ is illfounded, then $P_\psi(\gamma)$ is an illfounded proof tree only if $\xi <\gamma$ and $\gamma$ is admissible. Moreover, there is an infinite desending chain through $D_\psi(\xi)$ definable over $L_{\xi^+}$, so $P_\psi(\xi^{++})$ is an illfounded proof tree. 
Let $F_\psi$ be a pre-dilator which maps an ordinal $\gamma$ to the Kleene-Brouwer ordering of $P_\psi(\gamma)$. We again have $F_\psi(\gamma)$ is illfounded  only if $\xi < \gamma$ and $\gamma$ is admissible, an that $F_\psi(\xi^{++})$ is illfounded.
Since $T$ proves the existence of $\xi$ and $T$ extends $\Pi^1_1{-}\ca$, $T$ proves the existence of $\xi^+$, and so it also proves that $F_\psi$ is not a dilator. This shows that $\xi^+<\alpha$, as claimed.
\endproof
This completes the proof of the theorem.
\endproof

\begin{remark}
The implication from \eqref{TheoremCharacterizations122} to \eqref{TheoremCharacterizations121} in the proof of Theorem \ref{TheoremCharacterizations12} did not use the fact that $T$ extended $\Pi^1_1${-}$\ca$ in an essential way. In fact a similar argument shows the following variant for extensions of $\aca$ (or $\kp$): if $\alpha$ is admissible and is the least ordinal such that $L_\alpha\models S$ for some recursively enumerable set of $\Sigma_1$ sentences $S$ in the language of set theory, then $\alpha = s^1_2(T)$ for some $\Sigma^1_2$-sound, recursively enumerable extension of $\aca$.

We sketch the proof of this variant. Define $P^\psi = \{P^\psi(\gamma):\gamma\in\Ord\}$ as before and define $T$ as before, except that we add $\aca$ instead of $\Pi^1_1${-}$\ca$ as an axiom.
In order to show that $s^1_2(T)\leq\alpha$,
we use a Lemma \ref{LemmaSoundOmegaModel} to find a model $M$ of $T$ whose only ordinals are those ${<}\alpha$.
\end{remark}

\begin{theorem}
$\delta^1_2 = \sup\{s^1_2(T): T$ is a recursively enumerable, $\Sigma^1_2$-sound extension of $\aca\}$.
\end{theorem}
\proof
The upper bound is immediate from Lemma \ref{Lemmas12upperbound}; the lower bound is immediate from Theorem \ref{TheoremCharacterizations12}.
\endproof

\begin{proposition}\label{LemmaPi12Axiomatized}
Let $T$ be a $\Pi^1_2$-axiomatized, fully sound extension of $\aca$. Then, $s^1_2(T) = \omega_1^{ck}$.
\end{proposition}
\proof
By Lemma \ref{LemmaSoundOmegaModel} applied to $\omega_1^{ck}$ we can construct an $\omega$-model $M$ such that $M$ contains only recursive ordinals and satisfies all true parameter-free $\Pi^1_2$ sentences; in particular $M$ satisfies $T$. Thus, letting $T_M$ be the theory of $M$, we have $s^1_2(T)\leq s^1_2(T_M)$. Now, suppose that 
\[\text{$M \models $``$D$ is not a dilator,''}\]
so there is an ordinal $\alpha\in M$ such that $D(\alpha)$ is illfounded. Since $\alpha^+ = \omega_1^{ck}$ by choice of $M$, we must have $\alpha<\omega_1^{ck}$. This implies $s^1_2(T)\leq s^1_2(T_M)\leq\omega_1^{ck}$.

In order to complete the proof of the proposition, it suffices to show that $\omega_1^{ck} \leq s^1_2(\aca)$.
For this, suppose $\alpha<\omega_1^{ck}$. We have to find a pre-dilator $D$ such that $\aca\vdash$ ``$D$ is not a dilator,'' but the least counterexample to the illfoundedness of $D$ is strictly greater than $\alpha$.

Let $a$ be a recursive wellordering of $\mathbb{N}$ of length greater than $\alpha$. Let $D_{a\to \omega^*}$ be the pre-dilator from Lemma \ref{LemmaImplicationDilator}, where $\omega^*$ denotes an infinite descending chain. Thus, $D_{a\to \omega^*}$ is not a dilator and moreover $D_{a\to \omega^*}(x)$ is illfounded if and only if there is an embedding from $a$ to $x$.
Let $F$ be the pre-dilator given by
\[F(x) = C_a + D_{a\to\omega^*}(x),\]
where $C_a$ denotes the constant dilator with value $a$.
Then $F$ is not a dilator, and moreover $F(x)$  is illfounded if and only if there is an embedding from $a$ into $x$. It remains to show that 
\[\aca\vdash \exists x\, \lnot\wo(C_a + D_{a\to\omega^*}(x)).\]
We reason in $\aca$. There are two cases. If $a$ is illfounded, then for every ordinal $\gamma$, $a + D_{a\to\omega^*}(\gamma)$ is illfounded. If $a$ is wellfounded, then $D_{a\to\omega^*}(a)$ is illfounded, so $a + D_{a\to\omega^*}(a)$ is illfounded, as desired.
\endproof

The following proposition clarifies the inequality in Proposition \ref{Propo12Lesss12}.
\begin{proposition}
There are $\Sigma^1_2$-sound, recursively enumerable extensions $T, S$ of $\aca$ such that:
\begin{enumerate}
\item $s^1_2(T) = o^1_2(T)$;
\item $s^1_2(S) < o^1_2(S) < \infty$.
\end{enumerate}
\end{proposition}
\proof
For the first claim, we let $T$ be the theory obtained in the proof of Theorem \ref{TheoremSpectrumAdmissible} in the case $\alpha = \omega_1^{ck}$. This theory satisfies $o^1_2(T) = \omega_1^{ck}$ by the statement of the theorem. Moreover, the proof shows that the $\omega$-model $M$ 
obtained from applying Lemma \ref{LemmaSoundOmegaModel} to the case $\alpha = \omega_1^{ck}$
satisfies $T$. Since it
contains only ordinals ${<}\omega_1^{ck}$, we have $s^1_2(T) = \omega_1^{ck}$ by arguing as in Proposition \ref{LemmaPi12Axiomatized}.

For the second claim, we let $S$ be the theory obtained in the proof of Theorem \ref{TheoremSpectrumAdmissible} in the case $\alpha = \omega_2^{ck}$. Again by the theorem, we have $o^1_2(T) = \omega_2^{ck}$. This time, however, we let $M$ be the $\omega$-model $M$ 
obtained from applying Lemma \ref{LemmaSoundOmegaModel} \emph{to the case $\alpha = \omega_1^{ck}$}. 
Using the fact that $T$ is $\Pi^1_2$-axiomatized and that $M$ is $\Sigma^1_1$-correct with parameters by Lemma \ref{LemmaSoundOmegaModel}, we see that
$M \models S$, so again by an argument as in Proposition \ref{LemmaPi12Axiomatized} we conclude $s^1_2(T) = \omega_1^{ck}$.
\endproof

We conclude with the following result which concerns the $\Sigma^1_2$-ordinal of some notable $\Sigma^1_2$-sound theories. 
\begin{theorem}\ \label{TheoremS12Ordinals}
\begin{enumerate}
\item \label{TheoremS12Ordinals1} $s^1_2(\aca) = \omega_1^{ck}$;
\item \label{TheoremS12Ordinals1.5} $s^1_2(\kp) = \omega_1^{ck}$;
\item \label{TheoremS12Ordinals2} $s^1_2(\Pi^1_1$-CA$_0) = \omega_\omega^{ck}$;
\item \label{TheoremS12Ordinals3} $s^1_2(\Pi^1_2$-CA$_0) =$ least ordinal stable to the least non-projectible ordinal.
\end{enumerate}
\end{theorem}
\proof
Item \eqref{TheoremS12Ordinals1} is immediate from Proposition \ref{LemmaPi12Axiomatized}.
Item \eqref{TheoremS12Ordinals1.5} follows from  \eqref{TheoremS12Ordinals1}, since it implies that
\[\omega_1^{ck}\leq s^1_2(\aca) \leq  s^1_2(\kp) \leq \omega_1^{ck}.\]
Items \eqref{TheoremS12Ordinals2} and \eqref{TheoremS12Ordinals3}  both follow immediately from the proof of Theorem \ref{TheoremCharacterizations12}.
\endproof

\section{Further remarks}
Let us mention that the proof of Theorem \ref{TheoremPi12Dilator} adapts to prove the existence, uniqueness, and recursiveness of $|T|_{\Pi^1_n}$ up to bi-embeddability for $\Pi^1_n$-sound, recursively enumerable extensions of $\aca$. Here, $|T|_{\Pi^1_n}$ is defined naturally as the least $n$-ptyx into which all the $T$-provable $n$-ptykes are embeddable. How to derive from this an interesting classification theory for recursively enumerable extensions of $\aca$ according to their $\Pi^1_n$ consequences is not completely clear.\\

We have developed the abstract background for  proof-theoretic analyses of theories at the level of $\Sigma^1_2$ and $\Pi^1_2$. Although the approach via $\Sigma^1_2$ consequences might initially seem more natural, it appears that this type of analysis is not as informative as the approach via $\Pi^1_2$ consequences, as evidenced by the results of Section \ref{TheoremS12Ordinals}. 

The approach via $\Pi^1_2$ consequences, however, leads to the prospect of more informative and finer analyses, as well to a proof-theoretic classification theory which  provides insight not accessible purely on the basis of ordinal ($\Pi^1_1$) analyses. This is because ordinal analysis cannot distinguish between Categories B, C, and D at all. Theories in Category D satisfy $o^1_2(T) = \infty$ and thus are $\Pi^1_2$-sound. Hence it makes sense to carry out a full $\Pi^1_2$-analysis of such a $T$ and, in particular, to compute $|T|_{\Pi^1_2}$. A computation of the functor $|\aca|_{\Pi^1_2}$ is reported in the forthcoming article \cite{AgPa}.  

For theories in Categories B and C, one can carry out a \emph{quasi-$\Pi^1_2$} analysis, which consists of identifying the ordinals $o^1_2(T)$. This led to the Spectrum Problem, which is open for non-recursive inadmissible ordinals.
\begin{question}
Which, if any, non-recursive, inadmissible ordinals $\alpha<\delta^1_2$  are of the form $o^1_2(T)$ for some recursively enumerable extension $T$ of $\aca$?
\end{question}
Theories in Categories B and C are still well within the scope of study of Ordinal Analysis, and it makes sense to compare the results of Ordinal Analysis and the quasi-$\Pi^1_2$-analysis. By analyzing and slightly modifying the constructions of Theorems \ref{TheoremSpectrumRecursive} and \ref{TheoremSpectrumAdmissible}, one can have $|T|_{\Pi^1_1}$ and $o^1_2(T)$ vary arbitrarily, subject to the constraints of Theorem \ref{TheoremPi12Zero} and Lemma \ref{LemmaEpsilon}. It is worth asking if this holds in general. More precisely:
\begin{question}
Let $\gamma$ be recursive.
Suppose $o^1_2(T)$ is equal to some non-recursive, inadmissible ordinal $\alpha$. Is there a recursively enumerable extension $S$ of $\aca$ such that $o^1_2(S) = \alpha$ and $|S|_{\Pi^1_1} = \varepsilon_\gamma$?
\end{question}

We finish with a question on the  possible $\Sigma^1_2$-soundness ordinals of theories.
\begin{question}
Suppose $T$ is a $\Sigma^1_2$-sound, recursively enumerable extension of $\aca$. Must $s^1_2(T)$ be admissible?
\end{question}

\bibliographystyle{abbrv}
\bibliography{References}

\begin{thebibliography}{10}

\bibitem{Agetal}
J.~P. Aguilera, A.~Freund, M.~Rathjen, and A.~Weiermann.
\newblock {Ackermann and Goodstein Go Functorial}.
\newblock {\em Pacific J. Math.}
\newblock In press.

\bibitem{AgPa}
J.~P. Aguilera and F.~N. Pakhomov.
\newblock {A $\Pi^1_2$ Proof-Theoretic Analysis of Arithmetical Comprehension}.
\newblock Forthcoming.

\bibitem{Ba75}
J.~Barwise.
\newblock {\em Admissible Sets and Structures}.
\newblock Perspectives in Mathematical Logic. Springer-Verlag, Berlin, 1975.

\bibitem{Ca94}
J.~R.~G. Catlow.
\newblock A proof-theoretical analysis of ptykes.
\newblock {\em Arch. Math. Log.}, 33:57--79, 1994.

\bibitem{Fr19}
A.~Freund.
\newblock {$\Pi^1_1$-comprehension as a well-ordering principle}.
\newblock {\em Adv. Math.}, 355, 2019.

\bibitem{Fr71}
H.~M. Friedman.
\newblock Higher set theory and mathematical practice.
\newblock {\em Ann. Math. Logic}, 2(3):325 -- 357, 1971.

\bibitem{Ge36}
G.~Gentzen.
\newblock {Die Widerspruchsfreiheit der reinen Zahlentheorie}.
\newblock {\em Math. Ann.}, 112:493--565, 1936.

\bibitem{Gi81}
J.-Y. Girard.
\newblock {$\Pi^1_2$}-logic, part 1: Dilators.
\newblock {\em Ann. Math. Logic}, 21(2):75 -- 219, 1981.

\bibitem{GN85}
J.-Y. Girard and D.~Normann.
\newblock {Set Recursion and $\Pi^1_2$-Logic}.
\newblock {\em Ann. Pure Appl. Logic}, 28:255--286, 1985.

\bibitem{Goe31}
K.~Gödel.
\newblock {Über formal unentscheidbare Sätze der Principia Mathematica und
  verwandter Systeme I}.
\newblock {\em Monatsh. f. Math.}, 38:173--198, 1931.

\bibitem{Je03}
T.~J. Jech.
\newblock {\em Set Theory}.
\newblock Springer Monographs in Mathematics. Springer, 2003.

\bibitem{Je72}
R.~B. Jensen.
\newblock {The fine structure of the constructible hierarchy}.
\newblock {\em Ann. Math. Logic}, 4:229--308, 1972.

\bibitem{KW91}
A.~Kechris and W.~Woodin.
\newblock A strong boundedness theorem for dilators.
\newblock {\em Ann. Pure Appl. Logic}, 52(1):93 -- 97, 1991.

\bibitem{KS81}
J.~Ketonen and R.~Solovay.
\newblock {Rapidly Growing Ramsey Functions}.
\newblock {\em Ann. Math.}, 113(2):267--314, 1981.

\bibitem{KP82}
L.~Kirby and J.~Paris.
\newblock {Accessible Independence Results for Peano Arithmetic}.
\newblock {\em Bull. London Math. Soc.}, 14(4):285--293, 1980.

\bibitem{MaMo11}
A.~Marcone and A.~Montalb{\'a}n.
\newblock {The Veblen Function for Computability Theorists}.
\newblock {\em J. Symbolic Logic}, 76:576--602, 2011.

\bibitem{Ma75}
D.~A. Martin.
\newblock {Borel Determinacy}.
\newblock {\em Ann. Math.}, 102(2):363--371, 1975.

\bibitem{MoSh11}
A.~Montalb{\'a}n and R.~A. Shore.
\newblock {The Limits of Determinacy in Second-Order Arithmetic}.
\newblock {\em Proc. London Math. Soc.}, 104:223--252, 2011.

\bibitem{PW21}
F.~N. Pakhomov and J.~Walsh.
\newblock {Reducing $\omega$-model reflection to iterated syntactic
  reflection}.
\newblock To appear.

\bibitem{PH77}
J.~Paris and L.~Harrington.
\newblock {A Mathematical Incompleteness In Peano Arithmetic}.
\newblock In J.~Barwise, editor, {\em Handbook of mathematical logic}.
  North-Holland, 1977.

\bibitem{Po08}
W.~Pohlers.
\newblock {\em Proof Theory: The First-Step Into Impredicativity}.
\newblock 2008.

\bibitem{Ra99a}
M.~Rathjen.
\newblock {The Realm of Ordinal Analysis}.
\newblock In S.~B. Cooper and J.~Truss, editors, {\em Sets and Proofs}, pages
  219--279. Cambridge University Press, 1999.

\bibitem{RaWe11}
M.~Rathjen and A.~Weiermann.
\newblock {Reverse Mathematics and Well-Ordering Principles}.
\newblock In S.~B. Cooper and A.~Sorbi, editors, {\em Computability in Context:
  Computation and Logic in the Real World}. 2011.

\bibitem{Sa76}
G.~E. Sacks.
\newblock {Countable admissible sets and hyperdegrees}.
\newblock {\em Adv. Math.}, 20(2):212--262, 1976.

\bibitem{Sa90}
G.~E. Sacks.
\newblock {\em {Higher Recursion Theory}}.
\newblock Lecture Notes in Logic. Springer-Verlag, Berlin, New York, 1990.

\bibitem{Sch77}
K.~Sch{\"u}tte.
\newblock {\em Proof Theory}.
\newblock 1977.

\bibitem{Simpson}
S.~Simpson.
\newblock {\em {Subsystems of Second-Order Arithmetic}}.
\newblock 1999.

\bibitem{Sp55}
C.~Spector.
\newblock {Recursive Wellorderings}.
\newblock {\em J. Symbolic Logic}, 20:151--163, 1955.

\bibitem{Ta67}
G.~Takeuti.
\newblock {Consistency Proofs of Subsystems of Classical Analysis}.
\newblock {\em Ann. Math.}, 86(2):299--348, 1967.

\bibitem{Ta75}
G.~Takeuti.
\newblock {\em Proof Theory}.
\newblock 1975.

\end{thebibliography}

\end{document}